\documentclass[preprint]{elsarticle}

\usepackage{graphicx,epic,eepic}   
\usepackage{euscript}  
\usepackage{amsfonts}  
\usepackage{amssymb}
\usepackage{latexsym}
\usepackage{amscd}
\usepackage{multirow,setspace}
\newcommand{\field}[1]{\mathbb{#1}} 
 
\newtheorem{nthm}{Theorem}

\newtheorem{lem}[nthm]{Lemma} 
\newtheorem{prop}[nthm]{Proposition}

\newtheorem{defn}[nthm]{Definition}
 
\newtheorem{rmk}[nthm]{Remark}
\newtheorem{mainthm}{Theorem}

\def\sq#1#2{s_{#2}\left(#1\right)}

\def\eps{\epsilon}

\def\Tble{T_{b,\lambda,\eps}}
\def\oT{\cT_{\eps,\tau}}
 \def\Veut{V_{\eps,u,\tau}}

\def\mF{\mathfrak{F}}
\def\Um{\mathfrak{U}}
\def\mf{\mathfrak{f}}
\def\mc{\mathfrak{c}}
\def\um{\mathfrak{u}}

\def\mK{\mathfrak{K}}
\def\mL{\mathfrak{L}}

\def\mB{\mathfrak{B}}
\def\mm{\mathfrak{m}}
\def\mM{\mathfrak{M}}
\def\mT{\mathfrak{T}}

\def\mw{\mathfrak{w}}
\def\mW{\mathfrak{W}}
\def\mv{\mathfrak{v}}
\def\mV{\mathfrak{V}}
\def\my{\mathfrak{y}}
\def\mY{\mathfrak{Y}}
\def\mDu{\mathfrak{Du}}
\def\mDf{\mathfrak{Df}}

\def\mE{\mathfrak{E}}

\def\fC{\field{C}}

\def\fR{\field{R}}

\def\cA{\mathcal{A}}
\def\cB{\mathcal{B}}
\def\cC{\mathcal{C}}
\def\cD{\mathcal{D}}
\def\cE{\mathcal{E}}
\def\cF{\mathcal{F}}

\def\cH{\mathcal{H}}

\def\cL{\mathcal{L}}
\def\cN{\mathcal{N}}
\def\cM{\mathcal{M}}
\def\cT{\mathcal{T}}
\def\cR{\mathcal{R}}

\def\cP{\mathcal{P}}
\def\cS{\mathcal{S}}
\def\cU{\mathcal{U}}
\def\cW{\mathcal{W}}
\def\cV{\mathcal{V}}
\def\cX{\mathcal{X}}
\def\cZ{\mathcal{Z}}
\def\cO{\mathcal{O}}
\def\cQ{\mathcal{Q}}

\begin{document}
\large


\title{On Analytic Perturbations of a Family of Feigenbaum-like  Equations}


\author{Denis Gaidashev}

\address{
Department of Mathematics,
University of Uppsala, Uppsala,
 Sweden.\\
\ead{gaidash@math.uu.se}}



\begin{abstract}
We prove existence of solutions $(\phi,\lambda)$ of a family of of Feigenbaum-like equations 
\begin{equation}\label{family}
\phi(x)={1+\eps \over \lambda} \phi(\phi(\lambda x)) -\eps  x +\tau(x),
\end{equation}
where $\eps$ is a small real number and $\tau$ is analytic and small on some complex neighborhood of $(-1,1)$ and real-valued on $\fR$. The family $(\ref{family})$ appears in the context of period-doubling renormalization for area-preserving maps (cf. \cite{GK}).

Our proof is a development of ideas of H. Epstein (cf \cite{Eps1}, \cite{Eps2}, \cite{Eps3}) adopted to deal with some significant complications that arise from the presence of the terms $-\eps  x +\tau(x)$ in the equation $(\ref{family})$.  The method relies on a construction of novel {\it a-priori}  bounds for unimodal functions which turn out to be very tight. We also obtain good bounds on the scaling parameter $\lambda$. 

A byproduct of the method is a new proof of the existence of a Feigenbaum-Coullet-Tresser function. 

\end{abstract}

\maketitle


\setcounter{page}{1}

\section{Introduction} 

Since its original discovery  \cite{Fei1}, \cite{Fei2}, \cite{TC}, the Feigenbaum-Coullet-Tresser equation 
\begin{equation}\label{F_equation}
\phi(x)={1 \over \lambda} \phi(\phi(\lambda x)),
\end{equation}
whose solutions have attracted an extraordinary amount of interest. The study of this equation resulted in some spectacular breakthroughs in one-dimensional complex and real renormalization theory, which finally culminated in the proof of universality for unimodal maps in \cite{Lyu}.

In this paper we will consider the family of equations $(\ref{family})$, where $\epsilon \le 1$ and $\tau$ is small. This ``fixed point'' problem for the operator  
\begin{equation}\label{Ret}
\cR_{\eps,\tau}: \phi \mapsto {1+\eps \over \lambda} \phi \circ \phi \circ \lambda -\eps  \, id +\tau
\end{equation}
surfaces in the period doubling renormalization for two-dimensional maps. Specifically, we  have previously argued in \cite{GK}  that the area-preserving renormalization fixed point $F^*$ --- that is the area-preserving map that satisfies  $F^*=\Lambda_*^{-1} \circ F^* \circ F^* \circ \Lambda_*$, where $\Lambda_*$ is some coordinate change --- is almost one-dimensional in the sense that it is very close to the area-preserving H\'enon-like map
\begin{equation}\label{Henon_like}
H(x,y)=(\phi(x)-y,x-\phi(\phi(x)-y )),
\end{equation}
where $\phi$ is a solution of $(\ref{family})$ for $\epsilon=1$ and $\tau=0$.  An approach to an analytic proof of existence  of $F^*$ based on its proximity to the map $(\ref{Henon_like})$ has been also suggested in  \cite{GK}.  Proofs of existence of solutions of $(\ref{family})$ in this, interesting, case are, however, extremely technical. In this paper we concentrate on a simpler case of small $\epsilon$ and small $\tau$.

The problem $(\ref{family})$ will be reformulated and solved as a fixed point problem for an operator on some compact set of functions whose elements satisfy some  {\it a-priori} bounds.   A number of technical conditions in the proof  will be verified on a computer.  

The original computer-free proof of existence of the solution to the Feigenbaum-Coullet-Tresser equation $(\ref{F_equation})$ due to H. Epstein (cf \cite{Eps1}, \cite{Eps2}, \cite{Eps3}) was given  for  $\phi$'s that can be factorized as $\phi(x)=U(x^2)$, where $U$ is a  diffeomorphism. The presence of extra terms in the equation $(\ref{family})$ means that the solutions for $\epsilon \ne 0$ or $\tau \ne 0$ generally are not even functions anymore. We will, therefore, demonstrate existence of solutions on  the  {\it Epstein class} $\phi(x)=U(x)^2$ (see Definition \ref{Epstein_class} below).  Eventual existence of universal {\it a-priori} bounds  in the Epstein class ({\it ``beau''} bounds) is a seminal result of D. Sullivan \cite{Sul}.  We will, however, avoid a demonstration of existence of such bounds (in our case, for the operator $(\ref{Ret})$), by showing that there is a rather small compact and convex subset $\cA$ of function in the  Epstein class, which is {\it invariant} under the action of $\cR_{\eps,\tau}$.   

The {\it a-priori} bounds that we construct are new in the sense that they depend on the values of the derivative of the function at two points in the real slice of its domain as parameters; by doing this we were able to make the bounds very tight and significantly reduce the set $\cA$ which is guaranteed to contain the solution of $(\ref{family})$.

Another novelty of the proof is in the way we deal with complications that arise from the  presence of terms $\eps x$ and $\tau(x)$ in the equation $(\ref{family})$. The effect of these terms is a possible loss of univalence of $U^{-1}$. This in turn implies that one can not rely on {\it a-priori} bounds exclusively anymore, but rather one needs to make a set of assumptions on the derivative of $U^{-1}$, and show that these assumptions are, in a sense, reproduced.

As a bonus, the proof also demonstrates a certain property of stability of the space of solutions of $(\ref{family})$: for all sufficiently small $\epsilon$ and $\tau$ the diffeomorphic part $U$ of solutions  of $(\ref{family})$ on the  {\it Epstein class} $\phi(x)=U(x)^2$ lie in one and the same functional space,  independent of $\epsilon$ and $\tau$.

\section{Notation. Herglotz functions}

We will proceed with some definitions.

The upper and the lower half planes will be denoted as
$$\fC_\pm \equiv \{z \in \fC: \pm {\Im (z)} > 0 \}.$$

Let $J=(-l,r) \subset \fR$. Given such interval $J \subset \fR$, denote
$$
\fC(J) \equiv \fC_+ \cup \fC_- \cup J, \quad \fC_1 \equiv \fC((-1,1)).   
$$ 

Define $\cD_+(J,\theta)$ to be an open subset of $\fC_+$  bounded by a circular arc intersecting $\fR$ at the endpoints of $J$ at an angle $\theta$, and let $\cD_-(J,\theta)=\cD_+(J,\theta)^*$ where ${}^*$ stands for the complex conjugation.

Recall, that for every Riemann surface $U$, conformally isomorphic to the unit disk $D$, there exists a unique (upto a multiplication by a constant) metric, invariant under conformal automorphisms of $U$, called the {\it Poincar\'e metric}. We will denote the {\it Poincar\'e distance} in $U$ induced by this metric as ${\rm dist}_U$.  

In particular, the Poincar\'e metric on the unit disk $D$ is given by $d s =2 |d z| /( 1-|z|^2)$, and for any $z \in D$,
$${\rm dist}_D(0,z)=\log {1+|z| \over 1-|z|}.$$

The following Lemma is standard (see, for example, \cite{Eps2}).

\begin{lem}
The set 
$$\cD(J,\theta)=\cD_+(J,\theta) \cup \cD_-(J,\theta) \cup J.$$
is a {\it Poincar\'e neighborhood} of $J$ in $\fC(J)$, specifically
$$\cD(J,\theta)=\left\{z \in \fC:  {\rm dist}_{\fC(J)}(z,J) < \log{1+\tan(\theta/4) \over 1-\tan(\theta/4) }   \right\}.$$
\end{lem}

\noindent {\it Proof.}

The map $h \circ a$, where $a$ is the affine map of $J$ onto $(-1,1)$, and 
$$h(z)=\log{ 1+z \over 1-z},$$
is a conformal map of $\fC(J)$ onto
$$\cS_\pi=\left\{w \in \fC: |\Im (w)|<\pi \right\}.$$

$h \circ a$ maps $\cD(J,\theta)$ onto 
$$\cS_\theta=\left\{w \in \fC: |\Im (w)|<\theta \right\}.$$

Clearly, a point $w \in \cS_\theta$ iff ${\rm dist}_{\cS_\pi}(w,\fR) < {\rm dist}_{\cS_\pi}(i \theta,0)$. Now, map, $\cS_\pi$ onto the unit disk by
$$g(w)={\rm tanh}\left({w \over 4}  \right).$$

We get that $g \circ h \circ a$ maps $\fC(J)$ conformally onto the unit disk $D$, preserving the Poincar\'e distance between the points, therefore, for any $z \in \cD(J,\theta)$,
$${\rm dist}_{\fC(J)}(z,J) < {\rm dist}_D(0,{\rm tanh}(i \theta /4))={\rm dist}_D(0,\tan(\theta /4))=\log{ 1+\tan \theta/4 \over 1-\tan \theta /4}.$$
$\Box$

Given an interval $J \subset \fR$ and complex number $d$, $\Im(d)>0$, denote
$$
\fC(J,d) \equiv \fC(J) \setminus \{z \in \fC: \Re(z)=\Re(d), \Im(z) \ge \Im(d) \quad {\rm or} \quad \Im(z) \le -\Im(d)\}.
$$ 

$\fC(J,d)$ is a complex plane with four slits.

We will denote ${\cF}(\cD)$ the Banach space of functions holomorphic on a domain $\cD$ equipped with the uniform norm. A subset of functions in $\cF$ assuming their values in a set $\cE$, will be denoted by $\cO(\cD,\cE)$.

Suppose that $\cE$ is simply connected and open, $\cD$ and $\cE$ are real symmetric, the boundary of $\cE \cap \field{R}$ is $\{L,R\}$, and let ${\bf \mc}=\{\mc_1,\mc_2,\mc_3,\mc_4\}$ be a quadruple of real numbers, such that $\{\mc_1,\mc_2\} \in \cD$ and $\{\mc_3,\mc_4\} \in \cE$. We will further define
$${\cA}(\cD,\cE;{\bf \mc}) \equiv \left\{u \in {\cO}(\cD,\cE): u(z)=u(z^*)^*, u(\cD \cap \fC_\pm) \subset \overline{\cE \cap \fC_\pm}, u \left( \mc_1 \right)=\mc_3, u(\mc_2)=\mc_4 \right\}.
$$

It is a classical result that the set ${\cA}(\cD,\cE;{\bf \mc})$ is compact in $\cF(\cD)$ (cf \cite{Eps2}). Finally,
$$ \cA_1({\bf c}) \equiv \cA (\fC_1,\fC_1;{\bf c}),\quad \cA_{J,I,d}({\bf \mc}) \equiv \cA (\fC(J,d),\fC(I);{\bf \mc}),\quad \cA_{J,I,d,p}({\bf \mc}) \equiv \cA (\fC(J,d),\fC(I,p);{\bf \mc}).$$

 Clearly, a function $u$ in $\cA(\cD,\cE,{\bf \mc})$ is isomorphic to some $f \in \cA_1({\bf c})$ through unique conformal isomorphisms $\Phi_1$ and $\Phi_2$:
$$u=\Phi_2^{-1} \circ f \circ \Phi_1,$$ 
normalized so that 
$$\Phi_1(l)=-1, \quad \Phi_1(r)=1, \quad \Phi_1(a_1)=b_1 \quad {\rm and } \quad \Phi_2(L)=-1, \quad \Phi_2(R)=1, \quad \Phi_2(a_2)=b_2.$$
Here, $a_k$ and $b_k$ are some points that will be chosen conveniently, and 
$$c_{1,2}=\Phi_1(\mc_{1,2}), \quad c_{3,4}=\Phi_2(\mc_{3,4}).$$

 Functions in ${\cA}_1({\bf c})$, commonly referred to as Herglotz functions, admit the following integral representation:
\begin{equation}\label{int_rep}
f(z)-c_3= a (z-c_1) + \int d \nu(t) \left({1 \over t-z} - {1 \over t-c_1}   \right), 
\end{equation}
where $\nu$ is a measure supported in ${\fR} \setminus (-1,1)$. This integral representation can be used to obtain the following {\it a-priori} bounds  on ${\cA}_1({\bf c})$
\begin{eqnarray} 
\label{function_1}  {c_4-c_3  \over c_2-c_1} {1+c_2 \over 1+x} \ge &{f(x)-c_3 \over x-c_1} &\ge  {c_4-c_3  \over c_2-c_1} {1-c_2 \over 1-x}, \quad x \in (-1,c_2),\\
 \label{function_2}  {c_4-c_3  \over c_2-c_1} {1+c_2 \over 1+x} \le &{f(x)-c_3 \over x-c_1} &\le  {c_4-c_3  \over c_2-c_1} {1-c_2 \over 1-x}, \quad x \in (c_2,1),\\
\label{first_der} { 1+c_1 \over (x-c_1 )  (1+x) } \le &{f'(x) \over f(x)-c_3}& \le { 1-c_1 \over (x-c_1 )  (1-x) }  , \quad x \in (-1,1),\\
\label{second_der} {-2 f'(x) \over 1+x}  \le &f''(x)& \le {2 f'(x) \over 1-x},  \quad x \in (-1,1).
\end{eqnarray}

If $\Phi \arrowvert_\fR$ is a monotone function, then one can transfer the bounds $(\ref{function_1})$--$(\ref{second_der})$ to $\cA(\cD,\cE,{\bf \mc})$.

Finally, we will mention the following version of Schwarz Lemma  which will play an important role in our proofs below (cf \cite{Eps2}, \cite{Sul}, \cite{LY}):

\begin{lem}\label{Epstein_lemma}
Let $u:\fC_J \mapsto \fC_{J'}$ be a holomorphic  map such that $u(J) \subset J'$. Then for any $\theta \in (0,\pi)$,  $u(\cD_\pm(J,\theta)) \subset \cD_\pm(J',\theta)$.
\end{lem}

\section{Summary of main results}

We will now summarize the main findings of the paper in a somewhat abridged form.

\begin{mainthm}\label{main_thm_1}
Set 
\begin{eqnarray}
\label{I1}I_1&=&(-1.23,0.23), \quad \theta_1={4  \over 5} \pi,\\
\label{I2} I_2&=&(-1.63975634,1.63975634), \quad \theta_2=0.830267 \pi,\\
\label{I3}I_3&=&(-1.6760020,1.6760020), \quad \theta_3=0.830825 \pi,
\end{eqnarray}
and $\cD=\cD_1(I_1,\theta_1)$,  $\cE=\cD(I_2,\theta_2)  \cap \cD(I_3,\theta_3)$.

There are numbers $\delta>0$, $\varepsilon>0$, $\nu>0$ and $\rho$, such that for any $0 \le \eps \le \nu$, and any $\tau$ holomorphic on $\cE$, real-valued on $\fR$, and satisfying  $\sup_{z \in \cE} |\tau(z)| < \delta$,  $\sup_{z \in \cE} |\tau'(z)| < \varepsilon$, $\tau(0)=0$, there exists a function $\phi_{\eps,\tau}$, holomorphic on some complex neighborhood $\cO$ of $L=(-1,1)$, and satisfying $\phi_{\eps,\tau}(0)=1$, and a number $\lambda$, such that the following holds:

\begin{itemize}

\item[i)] $\phi_{\eps,\tau}$ and $\lambda$ solve the equation $(\ref{family})$ on ${\cO}$;

\item[ii)] $\phi_{\eps,\tau}$ has a unique quadratic critical point on $\cO$: $\phi_{\eps,\tau}(c+z)=O(z^2)$;

\item[iii)] the two inverse branches $\eta$ and $\zeta$ of $\phi_{\eps,\tau}$ can be factorized as
\begin{equation}\label{inverse_branches}
\eta(z)=u(T(-\sqrt{L(z)})), \quad \zeta(z)=u(T(\sqrt{L(z)})),
\end{equation}
$T$ and $L$ are  affine, and $ u$ belongs to a convex subset of ${\cA}(\cD,\cE;{\bf \mc}), \quad {\bf \mc}=\left(-1 / 2,0,0,1\right)$.


\end{itemize}
\end{mainthm}

Unsurprisingly, in the particular case of $\eps=0$, $\tau=0$, one can demonstrate that the factorized inverse of the Feigenbaum function $\phi^* \equiv \phi_{0,0}$ has much nicer analytic  properties. This is emphasized in our second result --- yet another proof of the existence of solutions of the  Feigenbaum-Coullet-Tresser equation --- which we state now in a simplified form:

\begin{mainthm}\label{main_thm_2}
Set  $d=0.5+0.352 i$, $p=0.69 i$ and
$$J=(-1.05,0.05), \quad I=(-1.1593855,1.1593855).$$

There exists a function $\phi^*$ analytic on some complex neighborhood  $\cE$ of $(-1,1)$ and satisfying $\phi^*(0)=1$, and a number $\lambda$, with the following properties:

\begin{itemize}

\item[i)] $\phi^*$ and $\lambda$ solve on $\cE$ the equation $(\ref{F_equation})$.

\item[ii)] $\phi^*$ has a unique quadratic critical point on ${\cE}$ at $0$: $\phi^*(z)=O(z^2)$;

\item[iii)] the two inverse branches $\psi$ and $\zeta$ of $\phi^*$ can be factorized as
\begin{equation}\label{inverse_branches_F}
\psi(z)=u(T(-\sqrt{1-z})), \quad \zeta(z)=u(T(\sqrt{1-z})),
\end{equation}
where $T$ is some explicit affine map, and $u$ belongs to a convex subset of $\cA_{J,I,d,p}({\bf \mc})$, ${\bf \mc}=(-1/2,0,0,1)$.

%

\item[iv)] $-0.40791 < \lambda < -0.38132$.

\end{itemize}

\end{mainthm}

We  emphasize that the proof supplies quite tight bounds on the scaling parameter $\lambda$. 

\section{Inverse branches. An operator on a compact space} \label{reduction}

In this section we will derive equations for the inverse branches of the solution of $(\ref{family})$.

We will look for this solution within a class of functions which are unimodal on some interval $I  \equiv [a,d] \supset [0,1]$, that is they have a unique critical point on $I$, and that this critical point $c$ is quadratic in the sense that  $\phi_{\eps,\tau}(x)=O((x-c)^2)$, and we will derive equations that the two inverse branches of such $\phi_{\eps,\tau}$ should satisfy.  Write 
$$\phi_{\eps,\tau}(x)=b-g(x-c), \quad b \equiv \phi_{\eps,\tau}(c),$$
then $(\ref{family})$ can be written as 
\begin{equation}\label{g_eq_2}
 g=F \circ g \circ \xi + \eps \!\! \quad \!\! id-\tau \circ (id +c),
\end{equation}
where
$$F(x)=b+c-{1+\eps \over \lambda} (b-g(b-c-x)), \quad  \xi(x)=\lambda x +c(\lambda-1).$$

We will now write a set of equations for the two inverse branches, $h$ and $f$, of $g$:
$$h: (0,g(d-c)) \mapsto (0,d-c), \quad f: (0,g(a-c)) \mapsto (a-c,0).$$ 

The inverse of  $(\ref{g_eq_2})$ on $(0,d-c)$ is the following set of equations for the inverse branches:
\begin{eqnarray}
\label{branch_1} f \circ F^{-1} \circ (id -\eps h+\tau  \circ (h+c))&=&\xi \circ h, \quad {\rm on} \quad (E,g(d-c)), \\
\label{branch_2} h \circ F^{-1} \circ (id -\eps h+\tau  \circ (h+c))&=&\xi \circ h, \quad {\rm on} \quad (0,E),
\end{eqnarray}
where $E \equiv g \left(c / \lambda-c\right)$.  The inverse of $(\ref{g_eq_2})$ on $(a-c,0)$ reads:
\begin{equation}\label{branch_3}
h \circ F^{-1} \circ (id - \eps f + \tau \circ (f+c))=\xi \circ f,  \quad {\rm on} \quad (0,g(a-c)).
\end{equation}

It is easy to check that, for example, functions $\phi_{\eps,0}$ for any  nonzero $\eps$ can not be even. We will, therefore, consider a larger class of functions (see, for example, \cite{EE}, \cite{Ya} in the context of critical circle maps) :

\begin{defn}\label{Epstein_class}
An orientation preserving interval homeomorphism $\phi: I \mapsto J$ belongs to  the {\it Epstein class}, if it extends to an analytic two-fold branched  covering of a topological disk $D \supset I$ onto the double-slit plane $\fC(J)$. A map $\phi$ in the Epstein class admits a factorization
$$\phi=q_c \circ U,$$
where $q_c(x)=x^2+c$, and $U$ is a univalent map of $D$ onto the complex plane with four slits which double covers $\fC(J)$ under the quadratic map $x \mapsto x^2+c$. 
\end{defn}

The Epstein class includes functions whose restriction to the real line is even:

\begin{lem}
Every $\phi$ which admits a decomposition $\phi=U \circ q_c$ on some topological disk $D \supset 0$,  with $U$ univalent on $q_c(D)$, is in the  Epstein class.
\end{lem}

From now on, we will consider functions in the Epstein class:
$$\phi(x)=U(x)^2,$$
and we will write
\begin{equation}\label{factorization}
{h=v \circ - \circ  s, \quad f=v \circ s}, 
\end{equation}
where $v$ is a diffeomorphism on $K \equiv (-\sqrt{g(d-c)},\sqrt{g(a-c)})$, $s(x) \equiv \sqrt{x}$ (the principle square root) and $-(x) \equiv -x$. A similar factorization has been used in \cite{Sul} and \cite{LY} to obtain {\it a-priori} bounds for a quadratic polynomial. With this factorization equations $(\ref{branch_1})$--$(\ref{branch_3})$ become
\begin{equation}
\label{linearizer} 
\xi \circ v = v \circ V, \quad 
V(x)= \left\{ 
-\sqrt{F^{-1} (x^2-\eps v(x)+\tau(v(x)+c))}, \quad  x \in [e,\sqrt{g(a-c)}),\atop
\phantom{--} \sqrt{F^{-1}  (x^2-\eps v(x)+\tau(v(x)+c))}, \quad  x  \in (-\sqrt{g(d-c)},e)
\right. .
\end{equation}

We will now formally introduce an operator which will be later shown to be defined on ${\cA}(\cD,\cE,{\bf \mc})$ for some choice of $\cD$, $\cE$ and  ${\bf \mc}=\left(-1/2, 0, 0,1 \right)$. The operator is defined through the following sequence of steps.

\begin{figure}[t]
 \begin{center}
\begin{tabular}{c c c}
$\!\!\!\!\!\!\!\!\!\!\!\!\!$ \resizebox{50mm}{!}{\includegraphics{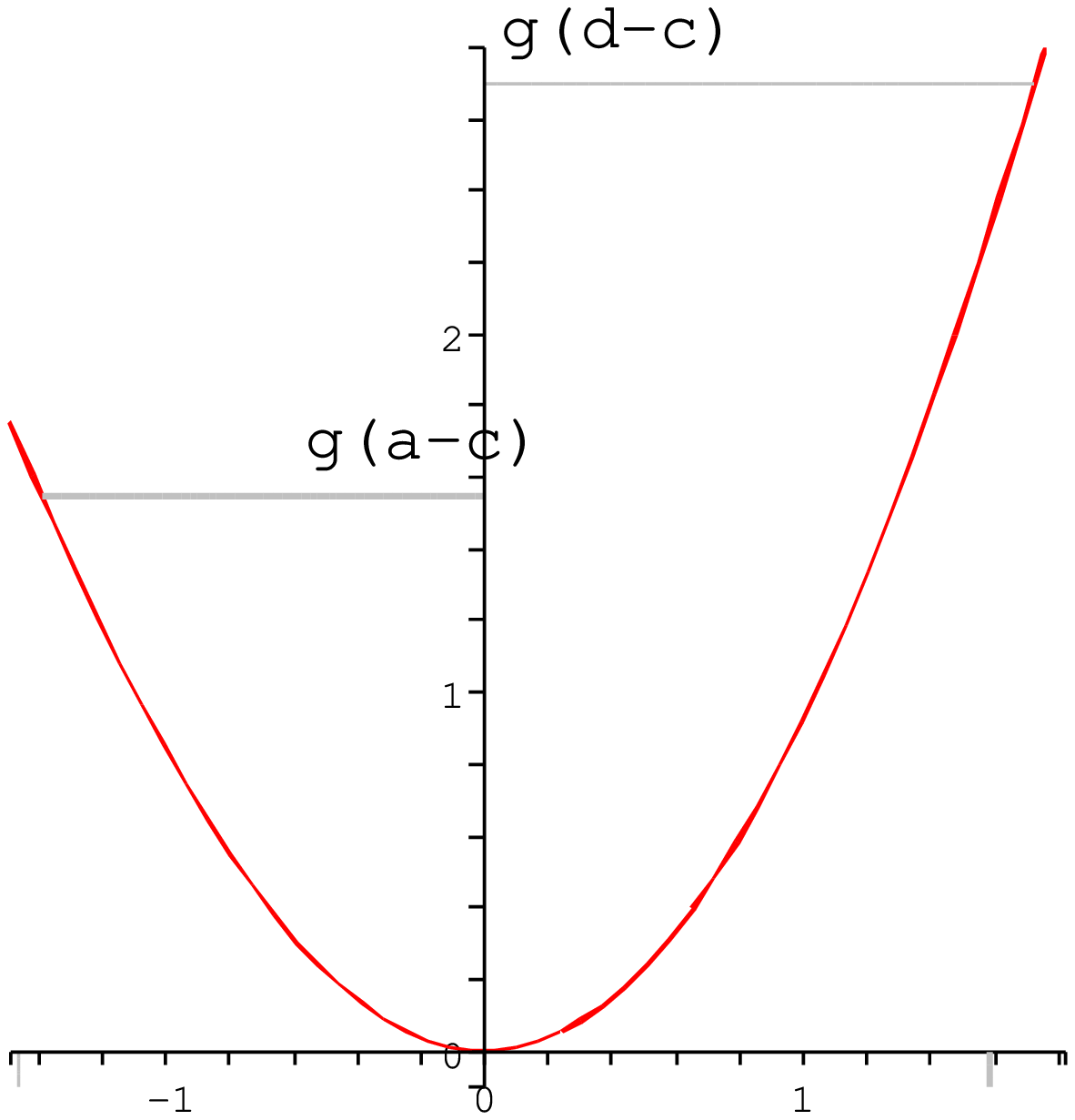}} &  $\!\!\!\!\!$ \resizebox{50mm}{!}{\includegraphics{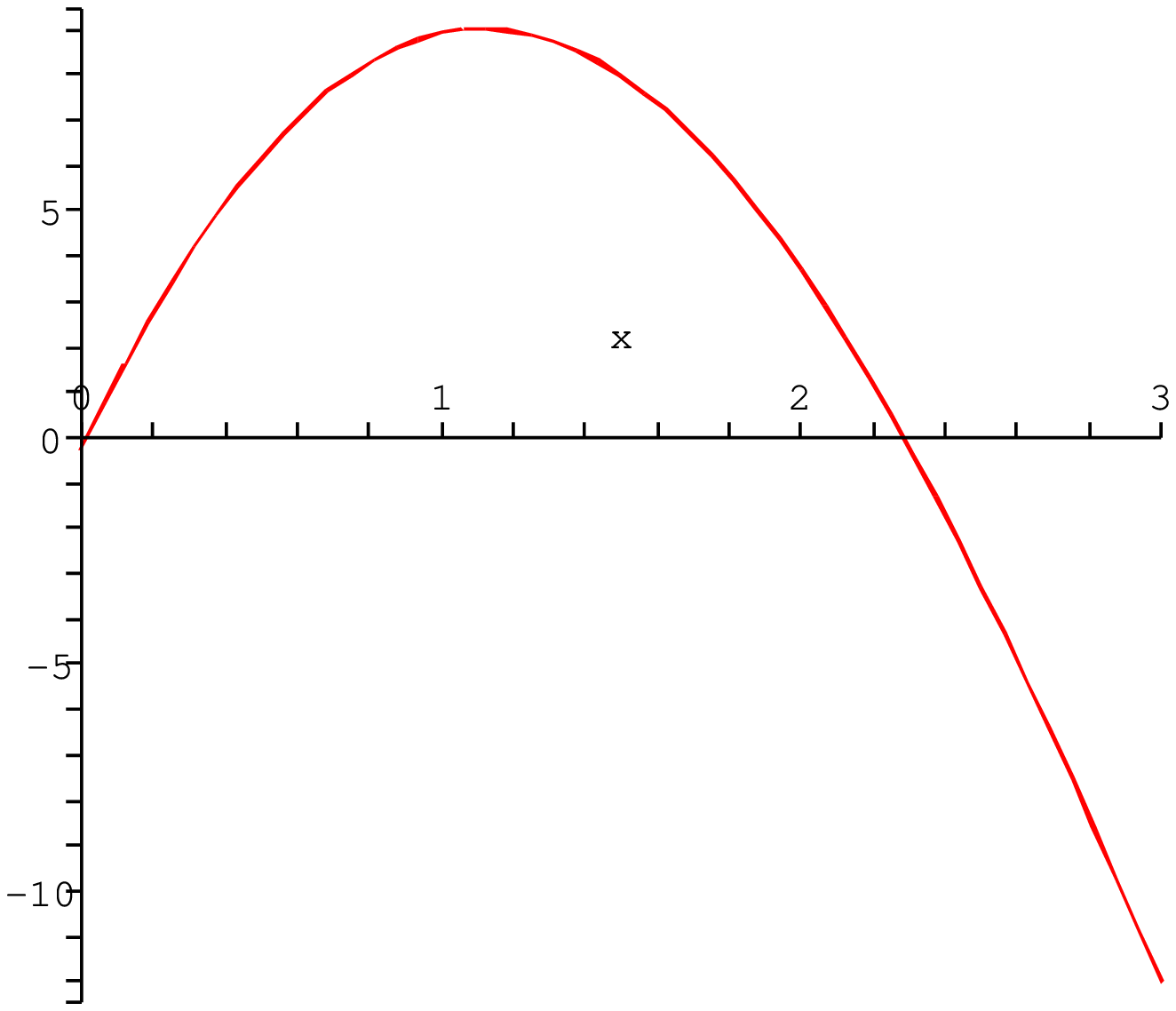}} & $\!\!\!\!\!$ \resizebox{50mm}{!}{\includegraphics{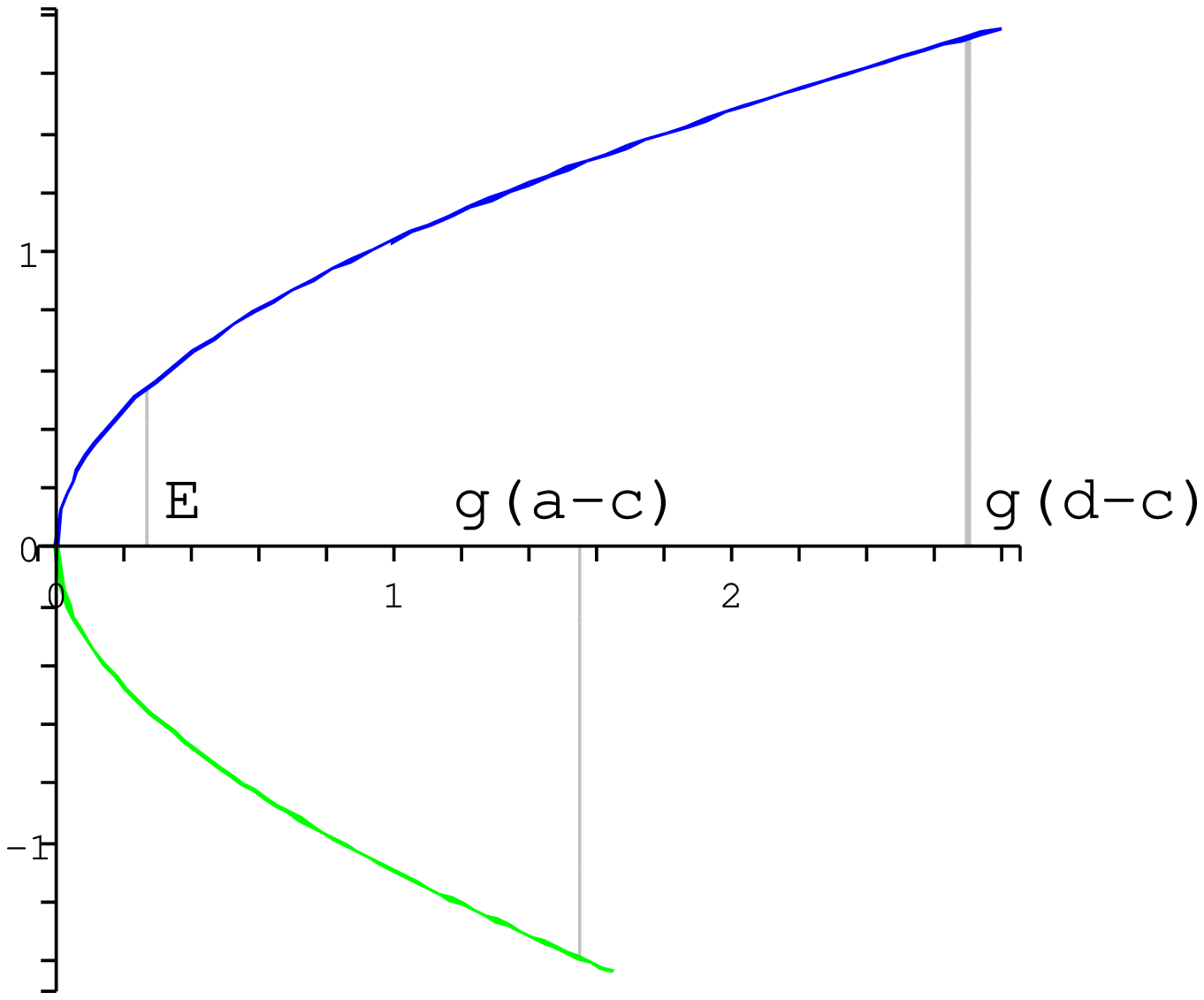}}\\
$\!\!\!\!\!\!\!\!\!\!\!\!\!$ a) &$\!\!\!\!\!$ b) & $\!\!\!\!\!$ c)
\end{tabular}
\caption{Function $g$ (in a)), function $F$ (in b)) and inverse branches $h$ and $f$ (in c)) for the solution $\phi$ of the equation $\phi(x)=2 \lambda^{-1}  \phi(\phi(\lambda x))-x$.} 
\end{center}
\end{figure}

\begin{itemize}
 \item[i)] Given $ u \in {\cA}(\cD,\cE,{\bf \mc})$, and a function $\tau$, holomorphic on $\cE \ni 0$, real-valued on $\fR$ and satisfying $\tau(0)=0$, find $b, \lambda$  and $e$ from the following set of equations  (for notational purposes,  we will use the symbol $s_b(x)$ for the function $\sqrt{b-x}$ through out the paper whenever convenient):
\begin{eqnarray}
\hspace{-1.0cm}\label{e_equation} -2 e  &=&  \alpha(b,\lambda,\eps) u'(\Tble(e))(\eps-\tau'(u(\Tble(e)))), \\
\hspace{-1.0cm}\label{l_equation} \lambda   &=&   u \left(\Tble \left(  \sq{u\left( \Tble  \left(-\sq{{ \eps  \lambda \over 1+\eps}+{\lambda^2 \over (1+\eps)^2}-{\lambda \over 1+\eps }\tau(1)}{b} \right) \right) }{b} \right) \right),\\
\hspace{-1.0cm}\label{b_equation} b  &=&  u \left( \Tble \left(  \sq{{ \lambda \over 1+\eps} ( b-e^2+\eps u(\Tble (e)) -\tau(u(\Tble(e)))}{b} \right) \right),
\end{eqnarray}
where $\alpha$, $\Tble$ and additional functions $\beta$ and $\gamma$ are given by
$$\alpha(b,\lambda,\eps)= { 1 \over 2 \beta(b,\lambda,\eps)-2 \gamma(b)}, \quad \beta(b,\lambda,\eps) =\sqrt{ b-{\lambda \over 1+\eps}}, \quad \gamma(b)=\sqrt{b-1},$$
$$\Tble(x) =-\alpha(b,\lambda,\eps) (x+\beta(b,\lambda,\eps)).$$
 
The affine transformation $T_{1,\lambda,0}$ will be also denoted by $T_\lambda$.

\item[ii)] Define for all $x \in T^{-1}_{b,\lambda,\eps}(\cD \cap \fR)$ 
\begin{equation}\label{V_equation}
 V_{\eps,u,\tau} (x) = {\rm sign}(e-x) \sq{{u  \left( \Tble  \left( - \left[ {w(\Tble(x))} \right]^{1 \over 2}\right)\right)}}{b},
\end{equation}
where
\begin{equation}\label{w_function}
 w (z) ={b- {\lambda \over 1+\eps} \left(b-\Tble^{-1}(x)^2+ \eps u(x) -\tau(u(x))\right)}.
\end{equation}

We will demonstrate that there is a choice of $\cD$ and $\cE$ such that  $V_{\eps,u,\tau}$ extends to a holomorphic function on $T^{-1}_{b,\lambda,\eps}(\cD)$.
 
\item[iii)] Set
\begin{equation} \label{T_op}
{\cT}_{\eps,\tau}[u](\Tble(z)) \equiv \lambda^{-1} u(\Tble(V_{\eps,u,\tau}(z))) .
\end{equation}

The operator $\cT_{0,0}$ will be denoted by $\cT$.

\end{itemize}

\begin{figure}[t]
 \begin{center}
\begin{tabular}{c c c}
 \resizebox{55mm}{!}{\includegraphics{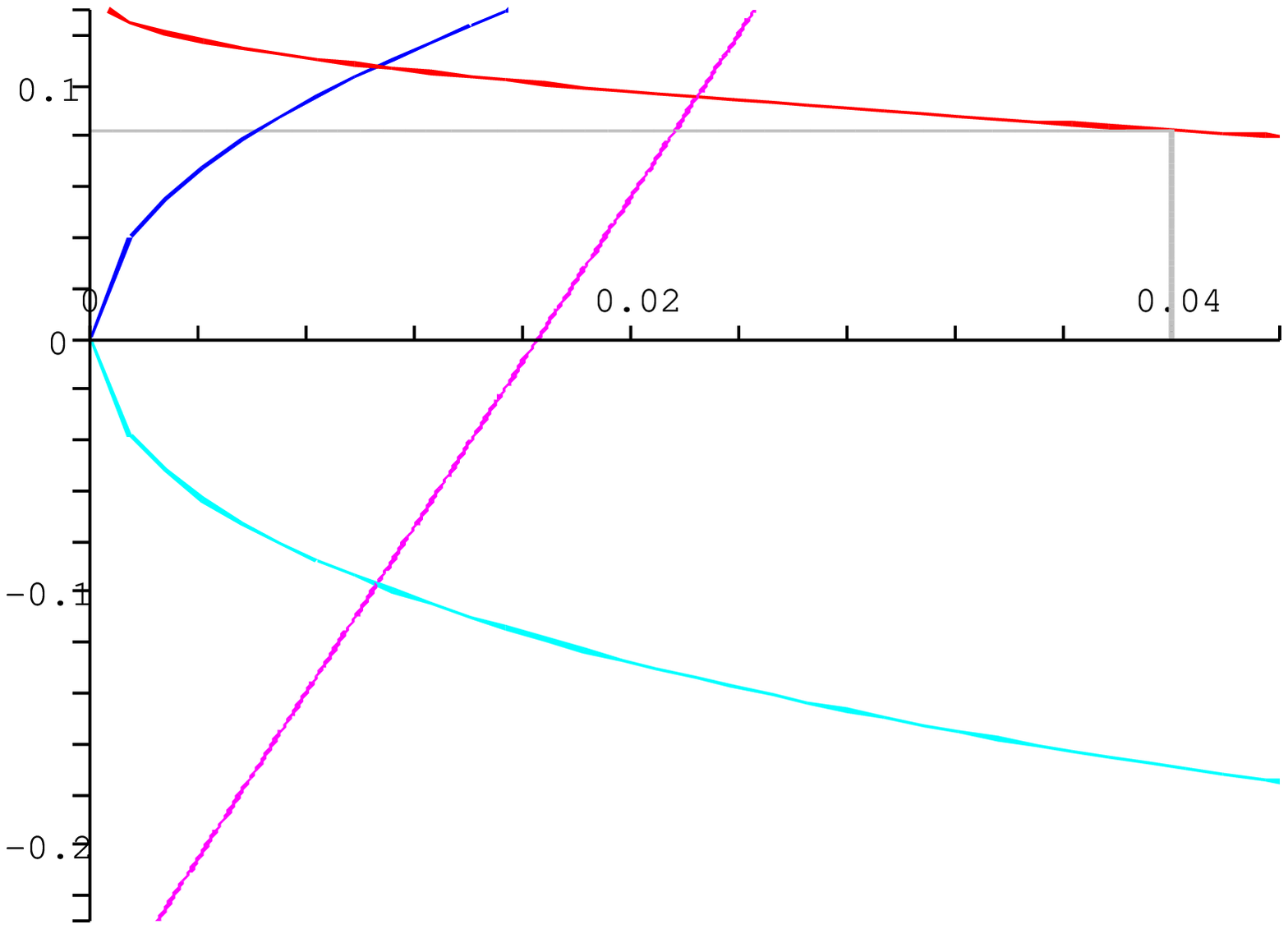}} & \phantom{aaaaa}  &  \resizebox{55mm}{!}{\includegraphics{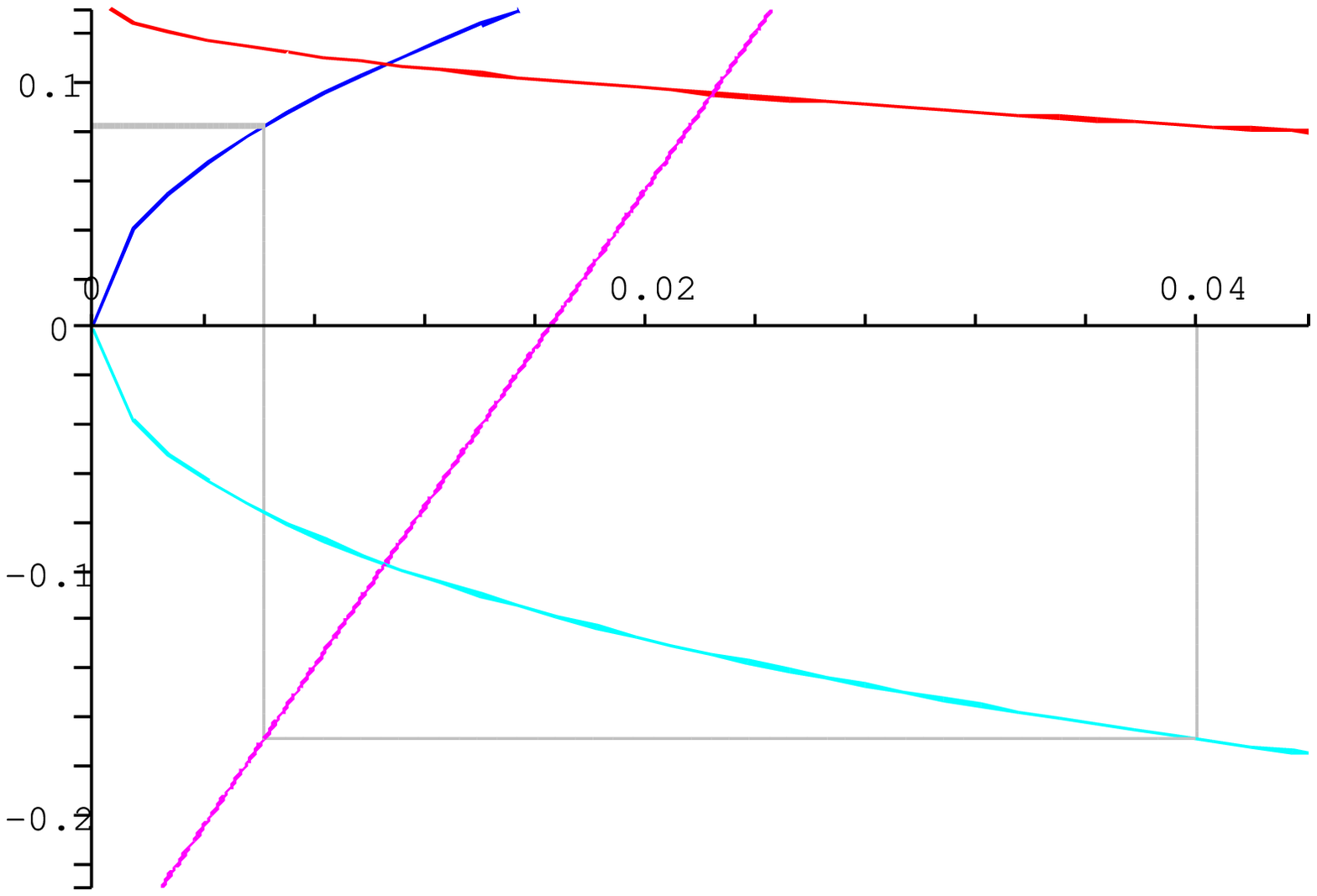}} \\
a) &\phantom{aaaaa} &  b)
\end{tabular}
\caption{An example of combinatorics in equalities $(\ref{branch_1})$ -$(\ref{branch_3})$ for a point in $(0,E)$ (equality $(\ref{branch_2})$): function $\xi \circ h$ is given in red, $id-h$ -- in cyan, $F$ -- in magenta, $h$ -- in blue; the image of the point under the right hand side of the equality is shown in a), under the left hand side -- in b).} 
\end{center}
\end{figure}

\begin{rmk}\label{RMK}
$ {} $

\noindent 1) {\rm Notice, that $\gamma=-\sqrt{b-1} \in (e,0)$ is the fixed point of $V_{\eps,u,\tau}$.}

\medskip

\noindent 2) {\rm  The normalization conditions $(\ref{e_equation})$--$(\ref{b_equation})$ ensure that $V_{\eps,u,\tau}$ is differentiable at $e$, and that } 
\begin{equation}
\nonumber  \cT_{\eps,\tau}[u](-1/2)=1, \quad \cT_{\eps,\tau}[u](0)=1.
\end{equation}

\medskip

\noindent 3) {\rm The  function $u$ is related to functions $v$, $\psi$, $h$ and $f$ appearing in the beginning of this Section  through the following equations:
\begin{eqnarray}
  \nonumber v(x)&=&u(-\alpha(x+\beta))-c, \\
  \nonumber h(x) &\equiv & =\psi(b-x)-c=u(\alpha (\sqrt{x}-\beta))-c, \quad x \in \left(0, \left[T^{-1}_{b,\lambda,\eps}(r)\right]^2\right), \\
  \nonumber f(x) &\equiv & u(\alpha (-\sqrt{x}-\beta))-c, \quad x \in  \left(0, \left[T^{-1}_{b,\lambda,\eps}(l)\right]^2\right).
\end{eqnarray}
}
\end{rmk}

We will show that for small $\eps$ and $\tau$, there is a choice of $\cD$ and $\cE$ such that that ${\cT}_{\eps,\tau}[u]$ is a continuous operator on $\cA(\cD,\cE,{\bf \mc}  )$. By compactness of the set $\cA(\cD,\cE,{\bf \mc}  )$ there is a function $u^*_{\eps,\tau} \in \cA(\cD,\cE,{\bf \mc}  )$ such that ${\cT}_{\eps,\tau}[u^*_{\eps,\tau}]=u^*_{\eps,\tau}$, which is equivalent to the set of equations $(\ref{branch_1})-(\ref{branch_3})$. In particular, $u^*_{\eps,\tau}$ is the ``factorized inverse'' (in the sense of Remark \ref{RMK} 3)) of a solution  of the equation $(\ref{family})$. 

\medskip

\noindent {\bf Remark 3.} Before we proceed with the proofs, we would like to emphasize two crucial difficulties that have forced us to modify the standard techniques that are commonly used to control inverse branches of unimodal maps (cf. \cite{Eps1}, \cite{Eps2}, \cite{Sul}, \cite{LY}).  

\noindent 1)  The terms $\eps x$ and $\tau(x)$ in the equation $(\ref{family})$ are responsible for the appearance of the terms $\eps u(\Tble(x))$ and $\tau(u(\Tble(x))))$ in $(\ref{V_equation})$--$(\ref{w_function})$. The effect of these terms is that one loses the benefit of estimating $u$, every time it enters the expression for $V_{\eps,u,\tau}$, only on a compact subset of its domain where one can use {\it a-priori} bounds. These terms do not appear in the Feigenbaum case ($\eps=\tau=0$) where this difficulty is absent. {\it In the case of nonzero $\eps$ and $\tau$ we are forced to  make assumptions on the range of $u$, and show that these assumptions are reproduced.}

\noindent 2) Another effect of terms $\eps u(\Tble(x))$ and $\tau(u(\Tble(x))))$ in $(\ref{V_equation})$ is that the derivative
$$\oT[u]'(z)=-\lambda^{-1} u'(\Tble(\Veut(\Tble^{-1}(z))) \alpha \Veut(\Tble^{-1}(z))'$$
can become zero since
$$\Veut(\Tble^{-1}(z))'=\ldots \times {1 \over \Veut(\Tble^{-1}(z))} \left({2 \over \alpha}  \Tble^{-1}(z)+\eps u'(z)-\tau'(u(z)) u'(z)\right)$$
can be zero. 

Notice, that $\Veut(\Tble(z))'$ is not zero at $\Tble^{-1}(e)$ (where the expression in parenthesis is equal to zero, cf $(\ref{e_equation})$): an application of the L'Hopital's rule shows that the derivative is finite at this point. However, it may be zero at other points on the real line where $2\alpha^{-1} \Tble^{-1}(z)+\eps u'(z)-\tau'(u(z)) u'(z)$ is zero. This would invalidate the argument since a function $\tilde{u} \equiv \oT[u]$ whose derivative is zero somewhere in the real slice of its domain generally is not in $\cA(\cD,\cE,{\bf \mc}  )$, in particular $\tilde{u}(\cD \cap \fC_\pm) \nsubseteq \overline{\cE \cap \fC_\pm}$. 

{\it We will deal with this problem by assuming an upper bound on the derivative $u'$ in the ``problematic'' subinterval of the real slice of $\cD$ so that $2 \alpha^{-1} \Tble^{-1}(z)+\eps u'(z)-\tau'(u(z)) u'(z)$ is guaranteed to be nonzero, and we will demonstrate that this bound is reproduced.}

\section{Yet another proof of existence of  the  Feigenbaum-Coullet-Tresser function}\label{FCC}

We will start by treating a simpler case of the Feigenbaum-Coullet-Tresser equation $(\ref{F_equation})$. The existence of solutions of the Feigenbaum-Coullet-Tresser equation is a well-established fact, and constitutes one of the most important results in one-dimensional renormalization theory. We will include this new proof here because  it illustrates some of the ideas used in a similar proof for equation $(\ref{family})$ in  the general case of nonzero $\eps$ and $\tau$  which could be otherwise obscured by technical details. 

Our proof follows the basic idea of H. Epstein of constructing an operator on a compact space of functions that admit {\it a-priori} bounds (cf \cite{Eps1}, \cite{Eps2}), but, at the same time, differs from it in that it is applicable to functions that are not necessarily even: $\phi(x)=U(x^2)$.

The case of the equation $(\ref{F_equation})$ is rather special. Suppose that $u \in \cA_{J,I,d,p}({\bf \mc})$ for some $J$, $I$, $d$  and $p$, and ${\bf \mc}=\left(-1 / 2, 0, 0,1 \right)$.  The set of  normalization conditions $(\ref{e_equation})$--$(\ref{b_equation})$ degenerates into simpler ones:
$$e =0, \quad b =u\left( T_{b,\lambda,0} \left[  -\sqrt{b- \lambda ( b-e^2)} \right] \right)= u\left( T_{b,\lambda,0} \left[  -\sqrt{b(1-\lambda)} \right] \right),$$
the last equation is clearly satisfied by $b=1$, since  $u\left( T_{1,\lambda,0} \left[  \sqrt{1-\lambda} \right] \right)=u(0)=1$. Then, the second normalization condition $(\ref{l_equation})$ becomes:

\begin{equation} \label{l_equation_F} 
\lambda =  u\left(T_\lambda \left( \sq{ u\left( T_\lambda \left[  -\sqrt{1-\lambda^2 } \right] \right) }{1}\right) \right), \quad {\rm where} \quad  T_\lambda \equiv T_{1,\lambda,0}.
\end{equation}

In the rest of this Section we will fix the following constants
$$l=1.05, \quad r=0.05, \quad m=1.1593855, \quad p=0.69 i,\quad d=0.5+0.3524 i,$$
and we will set $J=(-l,r)$, $I=(-m,m)$.  Furthermore, we will consider a smaller set of functions within $\cA_{J,I,d,p}({\bf \mc})$, specifically, functions that extend to $\fC(J,\bar{d}(t,s))$ with some $\bar{d}(t,s)\ge d$,  where 
$$s \equiv u'(0), \quad {\rm and} \quad  t \equiv u'(-1/2).$$

The set of such functions within  $\cA_{J,I,d,p}({\bf \mc})$ is clearly convex. We will refer to this set as $\cA_{J,I,\bar{d},p}({\bf \mc})$. The specific form of the continuous function $\bar{d}(t,s)$ will be described later.

The proof of the Proposition $\ref{central_prop_F}$ below is mildly computer assisted, and uses ``improved'' Herglotz bounds on $\cA_1({\bf c})$ transferred to ${\cA}_{J,I,\bar{d},p}({\bf \mc})$ with the help of the conformal isomorphisms 
$$\Phi_1(z,t,s)=A_1 { z- \Re{(\bar{d}(t,s))}  \over \sqrt{(z-\Re{(\bar{d}(t,s))})^2+ \Im{(\bar{d}(t,s))}^2}}+B_1,\quad \Phi_2(z)=A_2 { z- \Re{(p)}  \over \sqrt{(z-\Re{(p)})^2+ \Im{(p)}^2}}+B_2,$$
where $A_i$ and $B_i$ are found from the normalization conditions $\Phi_1(-l)=-1$, $\Phi_1(r)=1$,  $\Phi_2(0)=0$, $\Phi_2(1)=1$. $\Phi_1$ maps a plane with four slits $\fC(J,d)$ to a double slit plane $\fC_1$ conformally, while $\Phi_2$ is a conformal map of $\fC(I,p)$ to $\fC(I')$ for some interval $I'$.

The improvement of the Herglotz bounds (see Appendix A) uses the fact that $u'(0)$ and $u'(-1/2)$ can not be arbitrarily large, and that $u$ assumes its values in $\fC(I,p)$ (in particular, is bounded on $J$). We would like to point out that the derivatives $s=u'(0)$ and $t=u'(-1/2)$
play an important role as parameters in these new bounds. In particular, only a rather small region of the $(t,s)$-plane is admissible for $u$ such that $u(J) \subset I$. We will use that $\Phi_i \arrowvert_\fR$, $i=1,2$, are monotone, and will transfer the improved Herglotz bounds $\mf$ and $\mF$ (cf $(\ref{f_bounds})$) from $\cA_1({\bf c})$ to  ${\cA}_{J,I,\bar{d},p}({\bf \mc})$:
\begin{equation}\label{Ubounds}
\Um(x;t,s) \equiv \Theta_2(\mF\left(\Phi_1(x,t,s);t,s\right)), \quad \um(x;t,s) \equiv \Theta_2( \mf \left(\Phi_1(x,t,s);t,s\right) ),
\end{equation}
where $\Theta_2=\Phi_2^{-1}$. The next result is central to our proof.

\begin{prop}\label{central_prop_F}
Suppose $u \in {\cA}_{J,I,\bar{d},p}({\bf \mc})$. Then, there exists a bounded convex open set ${\cS}\subset \fR^2$, and two continuous functions $\cL_-(t,s)$ and  $\cL_+(t,s)$ such that the following holds whenever $(t,s) \equiv \left( u'\left(-{1 \over 2}\right),u'(0) \right) \in \cS$.

\begin{itemize}

\item[i)] There is a unique  $\lambda$,
\begin{equation}\label{lambda_b}
\cL_-(t,s) \le \lambda \le \cL_+(t,s),
\end{equation}
that solves $(\ref{l_equation_F})$. Furthermore, the map $u \mapsto \lambda$ is continuous.

\item[ii)]  The function $V_u \equiv V_{0,u,0}$ defined in $(\ref{V_equation})$ extends to a conformal map on ${\fC}(T^{-1}_\lambda(J),T^{-1}_\lambda(\bar{d}))$ that maps  ${\fC}(T^{-1}_\lambda(J),T^{-1}_\lambda(\bar{d})) \cap \fC_\pm$ into  ${\fC}(T^{-1}_\lambda(J),T^{-1}_\lambda(\bar{d})) \cap \fC_\mp$.

\item[iii)] Derivatives $\left(\cT[u]'(0), \cT[u]'(-1/2)\right)$ are also in $\cS$.

\end{itemize}
\end{prop}

\begin{figure}[t]
  \begin{center}
\resizebox{130mm}{!}{\includegraphics{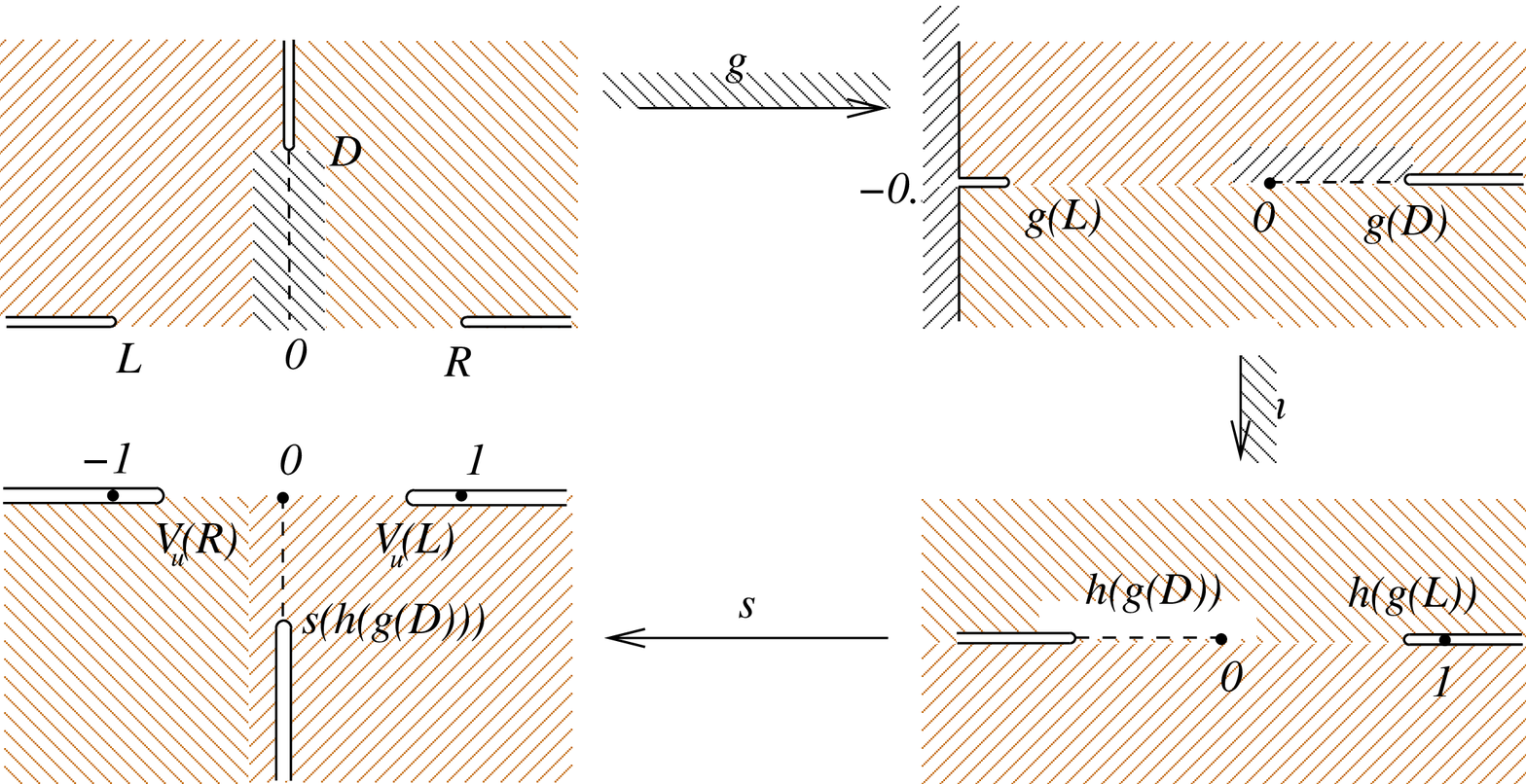}}
\caption{$V_u = s \circ h \circ g$ maps $\fC(K,D) \cap \fC_+$ into $\fC(K,D) \cap \fC_- $. The interval $H(t,s)$ and its images under maps $s$, $h$ and $g$ are given by dashed lines.}   
\label{fig2}
  \end{center}
\end{figure}


\noindent {\it Proof.}

\noindent i)  To demonstrate the claim of this part we consider the following function 
$$\mL(\lambda) \equiv \lambda \cT[u](0)=u\left(T_\lambda \left(\sq{u\left(T_\lambda \left(-\sq{\lambda^2}{1}  \right)\right)}{1} \right) \right),$$
and demonstrate that the function  $\lambda-\mL(\lambda)$ has exactly one zero in some interval $(\cL_-(t,s),\cL_+(t,s))$  for all $(t,s) \in \cS$. To this end we construct functions $\cL_\pm(t,s)$ so that the following holds
\begin{eqnarray}
\label{cL1} & &\cL_+(t,s)-\Um\left(T_{\cL_+(t,s)}\left(\sq{\Um\left(T_{\cL_+(t,s)}\left(-\sq{\cL_+^2(t,s)}{1}\right);t,s\right)}{1} ;t,s\right) \right) \ge 0, \\
\label{cL2} & &\cL_-(t,s)-\um\left(T_{\cL_-(t,s)} \left(\sq{\um\left(T_{\cL_-(t,s)} \left(-\sq{\cL_-^2(t,s)}{1}\right);t,s\right)}{1} ;t,s\right) \right)  \le 0,\\
& & 1-\mL'(\lambda) > 0 \quad {\rm for \quad all }\quad   \cL_-(t,s)  \le \lambda \le  \cL_+(t,s).
\label{d} 
\end{eqnarray}

The last inequality implies that if $\lambda$ is a zero of $1-\mL'(\lambda)$ then it is unique.

To demonstrate the inequalities we first choose a grid $\{ (t_i,s_k) \}$ of points on $\cS$, and  construct a set of numbers $\cL_\pm^{i,k}$ that satisfy $(\ref{cL1})$ and $(\ref{cL2})$ at $(t_i,s_k)$ numerically through a bisection procedure. We next define $\cL_\pm(t,s)$ over all of $\cS$ through a spline interpolation over points $\cL_\pm^{i,k}$. Finally, we verify that these functions $\cL_\pm(t,s)$ do satisfy $(\ref{cL1})$ and $(\ref{cL2})$  on all of $\cS$ using interval arithmetics.


\noindent ii)  Denote $D(t,s)=T^{-1}_{\lambda}(-\bar{d}(t,s))$, $L=T^{-1}_{\lambda}(r)$, $R=T^{-1}_{\lambda}(-l)$, $K=T^{-1}_{\lambda}(J)$ and  let $H(t,s)$ be the interval $(0, D(t,s))$ on the imaginary axis. 

First, we verify that $V_u$ is well-defined on $K$. For this, it is enough to check that 
$$1-\lambda (1-L^2)= 1-\lambda (1-R^2) >0,\quad T_{\lambda}\left(-\sq{\lambda(1-K^2)}{1} \right) \subset J,$$
$$0 < 1-u\left(  T_{\lambda}\left(-\sq{\lambda(1-L^2)}{1} \right)  \right) =  1-u\left(  T_{\lambda}\left(-\sq{\lambda(1-R^2)}{1} \right)  \right) < 1,$$
where the last inequality implies that $V_u(K) \Subset K$. These inequalities are verified on the computer for all $\cL_-(t,s) \le  \lambda \le \cL_+(t,s)$ and $(t,s) \in \cS$  using bounds $(\ref{Ubounds})$ (see \cite{Progs}).

Next, we shall extend $V_u$ first to  ${\fC}(K,D(t,s)) \cap \fC_+$ as $V_u = s \circ h \circ g$ where
$$g(z)=T_{\lambda}\left(-\sqrt{1-\lambda(1-z^2)} \right),\quad h(z)=1-u(z),\quad s(z)=- {\rm signum} (\Im(z)) \sqrt{z}$$
(see Fig. $\ref{fig2}$), and after that --- to  ${\fC}(K,D(t,s)) \cap \fC_-$ as $V_u(z) \equiv V_u^*(z^*)$ where $z^*$ signifies a complex conjugate of $z$. Functions $h$ and $g$ are not to be confused with those appearing in Section $\ref{reduction}$.

To this end, we first verify that $g$ maps ${\fC}(K,D(t,s) ) \cap \fC_+$ into the domain of $h$. For this we check that
$$0<1-\lambda (1-D(t,s)^2) <r \quad {\rm and} \quad  g(H(t,s)) \Subset  J,$$
for all $\cL_-(t,s) \le \lambda \le \cL_+(t,s)$, $(t,s) \in \cS$ and  $l$, $r$, $\bar{d}(t,s)$ and $p$ as in the condition.

 Next, notice that  $h \circ g$ maps quadrants $\fC_+ \cap \{z:\Re(z) \gtrless 0\}$ separately into $\fC_\pm$, which are mapped further  by  $s$ into $\fC_- \cap \{ z: \Re{z}  \lessgtr 0\}$. At the same time $h(g(H(t,s))) \subset \fR_-$, and therefore $h(g(H(t,s))$ is not in the domain of analyticity of $ \sqrt{}$. Therefore such $V_u$ is not defined on $H(t,s)$, but it is easily checked that it is continuous  across the interval $H(t,s)$; to be precise
$$\lim_{\epsilon \rightarrow 0} s(h(g(z-\epsilon)))=\lim_{\epsilon \rightarrow 0} s(h(g(z+\epsilon))), \quad z \in H(t,s),$$
and holomorphic in  $\fC_+ \cap \{z:\Re(z) \gtrless 0\}$. Therefore, by Morera's theorem, it is holomorphic in all of $\fC(K,D(t,s)) \cap \fC_+$.

To finish the verification of $V_u \left( \fC(K,Dt,s)) \cap \fC_+ \right) \subset \fC(K,D(t,s)) \cap \fC_- $ we have checked on the computer (see\cite{Progs})  that
\begin{equation}\label{ineq1}
 \lim_{\epsilon \rightarrow 0} |s(h(g(D(t,s)- \epsilon ))) \le \lim_{\epsilon \rightarrow 0}  |s(1-\um(g(D(t,s)\pm \epsilon );t,s))| \le |D(t,s)|
\end{equation} 
for all $\lambda$ as in $(\ref{lambda_b})$ and $(t,s) \in {\cS}$. 

Finally, 
$$V_u'(z)=-\lambda u'\left(T_{\lambda}(-\sq{\lambda(1-z^2)}{1})\right){ \alpha z  \over 2 V(z) \sq{\lambda(1-z^2)  }{1} }$$ 
and the only candidate for a zero of this derivative is $z=0$. However, $V_u(0)=0$, and an easy application of L'Hopital's rule demonstrates that $V'(0)=-\alpha\sqrt{u'(0) |\lambda|} \ne 0$. Therefore, $V_u$ is conformal.
 
\noindent iii) The proof of existence and invariance under $\cT$ of the set $\cS$ is practically identical to that of Lemma $\ref{lemma_S}$.

$\Box$

\medbreak

Existence of the fixed point of the the operator ${\cT}$ follows from the next

\begin{prop}\label{T0}
${\cT}$ is a continuous operator on ${\cA}_{J,I,\bar{d},p}({\bf \mc})$.
\end{prop}

\noindent {\it Proof.} Denote, as before, $s=u'(0)$ and $t=u'(-1/2)$, and let $P$ be the interval $(0,p)$ on the imaginary axis. To demonstrate that ${\cT}[u](J) \subset I$  whenever $u(J) \subset I$, and that ${\cT}[u](T_\lambda(H(t,s))) \subset P$ whenever  $u(T_\lambda(H(t,s))) \subset P$, it is enough to show that the functions
\begin{eqnarray}
{\cU}_1(\lambda,t,s) &\equiv &m-{1 \over \lambda} \Um\left(T_\lambda \left(\sq{\Um(T_\lambda(-\sq{\lambda^2}{1}) )}{1}\right) ;t,s\right), \\
{\cU}_2(\lambda,t,s) &\equiv &{1 \over \lambda} \um \left(T_\lambda \left(\sq{\um (T_\lambda(-\sq{\lambda^2}{1}) )}{1}\right) ;t,s\right)+m, \\
\label{q_ineq}\cQ(t,s) & \equiv & \Im{( p)} - \Im{ \left( \cT[u](\bar{d}(t,s)) \right) },      
\end{eqnarray}
are positive for all $(t,s) \in {\cS}$ and all $\lambda$ as in $(\ref{lambda_b})$.  The positivity of functions $\cU_i$ is verified on the computer (see \cite{Progs}).

To show that $\cQ(t,s) >0$, we use Lemma $\ref{Epstein_lemma}$. We first make some convenient choice of a Poincar\'e neighborhood $\cD_-((-x',x'),\theta')) \subset \fC(J,\bar{d})$  such that 
$$s(h(g(D(t,s)) \in  \cD_-((-x',x'),\theta')).$$
Then  Lemma $\ref{Epstein_lemma}$ guarantees that 
$$\lambda^{-1} u( T_\lambda (s(h(g(D(t,s)))))) \subset \lambda^{-1} \cD_-\left( (\um(-x';t,s),\Um(x';t,s)) ,\theta' \right) \equiv \cW(\lambda;t,s).$$
At this point we verify on the computer that the intersection of the set  $\cW(\lambda;t,s)$ with the imaginary axis is contained in the interval $(0,p)$ for all  $\cL_-(t,s) \le  \lambda \le \cL_+(t,s)$ and $(t,s) \in \cS$ --- that is we verify the inequality $(\ref{q_ineq})$. In fact, the function $\bar{d}(t,s)$ has been found, first, numerically, so that the inequality would be satisfied. This has been done over a grid of points in $\cS$ as a simultaneous bisection procedure for $\cL_\pm$ and $\bar{d}$ which finds some solutions of inequalities $(\ref{cL1})$, $(\ref{cL2})$ and $(\ref{q_ineq})$. Functions $\cL_\pm$ and $\bar{d}$ are next defined as some linear extrapolation over the grid, and the inequalities are checked again using the interval arithmetics.

 The above, together, with Proposition $\ref{central_prop_F}$ implies the claim.

$\Box$

\medbreak

 By the Tikhonov-Schauder theorem the operator $\cT$ has a fixed point in $\cA_{J,I,\bar{d},p}({\bf \mc})$. This completes the proof of Theorem $\ref{main_thm_2}$.

\section{General case $\eps, \tau \ne 0$.}\label{epsneq0}

In what follows, we will make the following choices: 
$$\cD=\cD(I_1,\theta_1), \quad \cE=\cD(I_2,\theta_2) \cap \cD(I_3,\theta_3),$$
where $I_1 \equiv (-l_k,r_k)$, and $\theta_k$ are as in $(\ref{I1})$--$(\ref{I3})$, and we will consider the corresponding space $\cA(\cD,\cE,{\bf \mc})$,
The point of considering an intersection of two Poincar\'e neighbourhoods as the target  set, rather than a single one, say $\cD(I_2,\theta_2)$, is that in our numerical experiments all choices of $\cD(I_1,\theta_1)$ and  $\cD(I_2,\theta_2)$, such that  the set $\cS$  of realizable $(t,s)$ is non-empty (and conveniently small) would lead to the target set $\cD(I_2,\theta_2)$ being too large for $\cT_{\eps,\tau}[u]$ to belong to the same space  $\cA(\cD(I_1,\theta_1),\cD(I_2,\theta_2),{\bf \mc})$. ``Clipping'' the target set by considering an appropriate intersection of two Poincar\'e neighbourhoods has enabled us to demonstrate the invariance of the space  $\cA(\cD,\cE,{\bf \mc})$ under $\cT_{\eps,\tau}$.

The double slit plane $\fC_1$ is isomorphic to Poincar\'e neighbourhoods $\cD(I_k,\theta_k)$ via  conformal isomorphisms
\begin{equation}\label{theta_k}
\Theta_k \equiv q_k \circ \sigma_k \circ m_k \circ  \zeta,
\end{equation}
where
\begin{eqnarray}
\nonumber \zeta (z) &\equiv&  { \sqrt{1+z}-\sqrt{1-z} \over \sqrt{1+z} + \sqrt{1-z} }, \hspace{12mm} m_k (z) \equiv  { z+a_k \over 1 - a_k z}, \\
\nonumber \sigma_k (z)  &\equiv& { (1+z)^{\kappa_k}-(1-z)^{\kappa_k}  \over (1+z)^{\kappa_k} + (1-z)^{\kappa_k} }, \hspace{7mm}  q_k (z) \equiv {l_k+r_k \over 2} z +{r_k-l_k \over 2},
\end{eqnarray}
where $-l_k$ and $r_k$ are the left and the right end points of intervals $I_k$,  and $ \kappa_k \equiv 2- {2 \theta_k / \pi}$.

With a little bit of  work, one can check that the transformation $\zeta$ maps $\fC_1$ onto the unit disk, $m_k$ is the normalizing Moebius transformation, $\sigma_k$ maps the unit disk onto $\cD((-1,1),\theta_k)$, and, finally, $q_k$ maps  $\cD((-1,1),\theta_k)$ onto  $\cD(I_k,\theta_k)$. Constants $a_k$ in the normalizing Moebius transformations $m_k$ are defined through the conditions  $\Theta_1(0) =-1/2$, $\Theta_2(0)=\Theta_3(0)=0$.



A function  $u$ in $\cA(\cD,\cE,{\bf \mc})$ can be now  factorized  as 
$$u= \Theta_2 \circ  f_2 \circ \Phi_1= \Theta_3 \circ f_3 \circ \Phi_1,$$
where 
$$f_k \in \cA_1({\bf c^k}), \quad  {\bf c^k}=(\Phi_1(\mc_1),\Phi_1(\mc_2),\Phi_k(\mc_3),\Phi_k(\mc_4)), \quad k=2,3.$$

 Therefore, according to Schwarz Lemma $\ref{Epstein_lemma}$, if $f_k \in \cA_1({\bf c^k})$ and an interval $J$ are such that $f_k(J) \subset J'_{k}$ then 
$$u( \Theta_1(\cD(J,\theta))) \subset \Theta_2(\cD(J'_2,\theta)) \cap  \Theta_3(\cD(J'_3,\theta)).$$

Furthermore, one can use the fact that $\Theta_k \arrowvert_\fR$ are monotone functions to transfer the improved Herglotz bounds $(\ref{f_bounds})$ from $\cA_1({\bf c^k})$ to $\cA(\cD,\cE,{\bf \mc})$:
\begin{eqnarray}
\label{U_bounds}\Um(x;t,s) &\equiv& \min \left\{ \Theta_2\left(\mF_2\left(\Phi_1(x);t,s\right)\right),\Theta_3\left(\mF_3\left(\Phi_1(x);t,s\right)\right) \right\}, \\
\label{u_bounds}\um(x;t,s) &\equiv& \max \left\{ \Theta_2\left(\mf_2\left(\Phi_1(x);t,s\right)\right),\Theta_3\left(\mf_3\left(\Phi_1(x);t,s\right)\right) \right\}.
\end{eqnarray}

We have implemented bounds $(\ref{U_bounds})$--$(\ref{u_bounds})$ on the computer, and used them in our proofs.

As in the previous section, we will consider a subset of $\cA(\cD,\cE,{\bf \mc})$ by allowing the real slices of the target sets $\cD(I_2,\theta_2)$ and $\cD(I_3,\theta_3)$ to be functions of $u'\left(-{1 \over 2} \right),u'(0))$:
\begin{eqnarray}
\nonumber I_2(t)&=&\left(\phantom{I_1^1}\!\!\!\!0.16(t-t^*)(0.5-l_1)-m_2,0.16(t-t^*)(r_1+0.5) +m_2 \right),\\
\nonumber I_3(t)&=&\left(\phantom{I_1^!}\!\!\!\!3.5(t-t^*)(0.5-l_1)-m_3,3.5(t-t^*)(r_1+0.5) +m_3\right),
\end{eqnarray}
where $m_2= 1.63825$, $m_3= 1.6430509$, and $t^*=1.9142899327$, $s^*=2.2366548836$ are approximate values of the derivatives  $u'\left(-{1 \over 2} \right)$ and $u'(0))$ for the fixed point of the operator $\cT$ computed numerically.

The subset $\{u \in \cA(\cD,\cE,{\bf \mc}): u(\cD) \subset  \cD(I_2(t),\theta_2) \cap \cD(I_3(t),\theta_3)\}$ is convex: if $u_1$ and $u_2$ are any two such functions and $(t_1,s_1)=\left( u_1'\left(-1 / 2 \right),u_1'(0)) \right)$ and $(t_2,s_2)=\left( u_2'\left(-1 / 2 \right),u_2'(0)) \right)$ are their derivatives, then any function $p u_1+(1-p) u_2$, $p \in (0,1)$, is also in the same subset. Indeed, if $z_1 \in \cD(I_k(t),\theta_k)$, $k=2,3$, and  $z_2 \in \cD(I_k(t),\theta_k)$, $k=2,3$, then  elementary geometric considerations demonstrate that $p z_1 +(1-p) z_2 \in \cD(I_k(p t_1+(1-p) t_2),\theta_k)$ (the fact that $|I_k(t)|$ is constant and independent of $t$ is important here).


We shall now proceed to describe a set $\tilde{\cS}$ of realizable derivatives $(u'\left(-{1 \over 2} \right),u'(0))$: 

\begin{lem}\label{lemma_S}
Suppose that $u \in {\cA}(\cD,\cE,{\bf \mc})$, and, furthermore, 
$$u(\cD) \subset \cD(I_2(t,s),\theta_2) \cap \cD(I_3(t,s), \theta_3).$$

Then there are four curves $(t,\cZ_2(t))$, $(t,\cZ_3(t))$, $(t,\cC_2(t)$ and  $t=t^*-0.0004$  in the $(t,s)$-plane that bound a convex open set $\tilde{{\cS}}$, such that
$$\left(\cT[u]'\left(-{1 \over 2} \right),\cT[u]'(0)\right) \subset \tilde{\cS}, \quad {\rm whenever} \quad \left(u'\left(-{1 \over 2} \right),u'(0)\right) \in \tilde{\cS}.$$
\end{lem}

\noindent {\it Proof.} See the Appendix B for the proof. 

$\Box$

The following Proposition shows that the space  $u \in {\cA}(\cD,\cE,{\bf \mc})$ is invariant under $\cT_{\eps,\tau}$.

\begin{prop}\label{central_prop}

There exist $\delta>0$, $\varepsilon>0$, $\nu>0$, and $C>0$ and $\sigma>0$, satisfying $C > \sigma \nu^2$, such that whenever 
\begin{itemize}
\item[1)] $u \in {\cA}(\cD,\cE,{\bf \mc})$, $u(\cD) \subset \cD(I_2(t,s),\theta_2) \cap \cD(I_3(t,s), \theta_3)$;
\item[2)] $\tau$ is a holomorphic function on $\cE$, real-valued on $\fR$, satisfying
$$\tau(0)=0, \quad  \sup_{z \in \cE} |\tau(z)| \le \delta, \quad  \sup_{z \in \cE} |\tau'(z)| \le \varepsilon;$$
\item[3)]for all $x \in (0,r_1)$ \begin{equation}\label{der_bound}
u'(x) \le  \omega + \rho x, \quad {\rm where} \quad  \omega=13, \quad \rho=20;
\end{equation}

\item[4)] $\eps$, the parameter in the operator $\cT_{\eps,\tau}$, is less than $\nu$;
\end{itemize}
there are two piecewise linear function $\cL_-(t,s)$ and $\cL_+(t,s)$,
, and two constant $\cB_- \equiv 1+\sigma \eps^4$ and $\cB_+ \equiv  1+C \eps^2$,  such that the following  holds
\begin{itemize}

\item[i)] there is a triple $(e,b,\lambda)$ that solves equations $(\ref{e_equation})$--$(\ref{l_equation})$, and satisfies
\begin{eqnarray}
\label{e_bounds} -\gamma(b) \ge  &e& \ge -\beta(b,\lambda,\eps), \\
\label{l_bounds}  \cL_+(t,s) \ge  &\lambda& \ge \cL_-(t,s),\\
\label{b_bounds}  \cB_+ \ge & b &\ge \cB_-,
\end{eqnarray}
where $t=u' \left( -1 / 2 \right)$, $s=u' \left( 0 \right)$. Furthermore, the map $u \mapsto  (e,b,\lambda)$ is continuous, while the solution $e$ of $(\ref{e_equation})$ is unique;

\item[ii)] $\oT[u]'$ also admits the bound $(\ref{der_bound})$;

\item[iii)]  the function $V_{\eps,u,\tau}$ extends to a holomorphic function on $T^{-1}_{b,\lambda}(\cD)$ that maps $T^{-1}_{b,\lambda}(\cD) \cup \fC_\pm $ compactly into  $T^{-1}_{b,\lambda}(\cD) \cup \fC_\mp $;

\item[iv)] $\oT[u] \in \cA(\cD,\cE, {\bf \mc})$, $u(\cD) \subset \cD(I_2(t,s),\theta_2) \cap \cD(I_3(t,s), \theta_3)$.

\end{itemize}
\end{prop}

We do not prove uniqueness of the solution $(b,\lambda)$, although this seems possible (with significantly more effort). We conclude that
$$\eta(z)=u_\infty(T_{b_\infty,\lambda_\infty}(-\sqrt{b_{\infty}-z})), \quad \zeta(z)=u_\infty(T_{b_\infty,\lambda_\infty}(\sqrt{b_{\infty}-z})),$$
are the factorized inverses of a solution $\phi_{\eps,\tau}$ of $(\ref{family})$ on some complex neighborhood of
$$u_\infty(I_1) \supset \left( \max_{(t,s) \in \tilde{\cS}} \Um(-l_1;t,s),\min_{(t,s) \in \tilde{\cS}} \um(r_1;t,s) \right) \supset (-1,1).$$


\subsection{Proof of part i) of Proposition  $\ref{central_prop}$.}\label{part1}

To demonstrate $(\ref{e_bounds})$ we introduce a function
$$ \mE(x;\lambda,b,\eps,\tau)  \equiv - u'(\Tble(x)) {\alpha(b,\lambda,\eps) \over 2} \left(\eps-\tau'\left( u(\Tble(x)) \right) \right).$$
Notice, that
\begin{eqnarray}
\nonumber \mE(-\gamma(b);\lambda,b,\eps,0)  &\equiv& -\eps u'(\Tble(-\gamma(b))) {\alpha(b,\lambda,\eps) \over 2}  =  - \eps t {\alpha(b,\lambda,\eps) \over 2}, \\
\nonumber \mE(-\beta(b,\lambda,\eps); \lambda,b,\eps,0) & \equiv & -\eps u'(\Tble(-\beta(b,\lambda,\eps))) {\alpha(b,\lambda,\eps) \over 2}  =  -\eps s {\alpha(b,\lambda,\eps) \over 2}.  
\end{eqnarray}

Since
$$\left(-\beta(b,\lambda,\eps),-\gamma(b) \right) \supset \left( -\sqrt{\cB_--\lambda}, -\sqrt{\cB_+-1} \right)=\left(-\sqrt{1+\sigma \eps^4 +|\lambda|}, -C^{1 \over 2} \eps \right),$$
for sufficiently small $\eps$ and for 
\begin{equation} \label{C1}
C  < \left(t {\alpha(b,\lambda,\eps) \over 2} \right)^2,
\end{equation}
the following holds
$$  - \eps t {\alpha(b,\lambda,\eps) \over 2}  < - C^{1 \over 2} \eps, \quad {\rm and} \quad -\sqrt{1+\sigma \eps^4+|\lambda|} < -\eps s {\alpha(b,\lambda,\eps) \over 2},$$
and the interval 
\begin{equation}
\label{e_cont}\left( \mE(-\beta(b,\lambda,\eps);\lambda,b,\eps,0) ,\mE(-\gamma(b);\lambda,b,\eps,0)   \right) \Subset(-\beta(b,\lambda,\eps),-\gamma(b)).
\end{equation}

 Since $\mE$ is clearly continuous in $\tau'$ at $\tau'=0$, there is a $\varepsilon>0$ such that the same containment $(\ref{e_cont})$ holds for all $\tau$ that satisfy  $\sup_{z \in \cE} |\tau'(z)| \le \varepsilon$.

 To show $(\ref{l_bounds})$--$(\ref{b_bounds})$ we consider two functions 
\begin{eqnarray}
\nonumber \mL_{u,\tau}(\lambda,b;\eps) &\equiv& u\left(\Tble\left(\sq{u\left(\Tble\left(-\sq{{\eps \lambda \over 1+\eps} +{\lambda^2 \over (1+\eps)^2}-{\lambda \over 1+\eps} \tau(1)}{b}  \right)\right)}{b} \right) \right), \\
\nonumber  \mB_{u,\tau}(\lambda,b;e,\eps) &\equiv& u \left(\Tble \left(-\sq{{\lambda \over 1+\eps} \left(b-e^2+\eps u(\Tble(e))+ \tau(u(\Tble(e))) \right) }{b} \right) \right),
\end{eqnarray}
and demonstrate that the map $(\lambda,b) \mapsto \left(\mL_{u,\tau}(\lambda,b;\eps),\mB_{u,\tau}(\lambda,b;e,\eps) \right)$ maps  the parallelogram $(\ref{l_bounds})$--$(\ref{b_bounds})$ in the $(\lambda,b)$-plane into itself for all $(t,s) \in \tilde{\cS}$, all $e$ as in $(\ref{e_bounds})$ and all $u \in {\cA}(\cD,\cE,{\bf \mc})$.  To this end, we first show that $\cL_+(t,s)-\mL_{u,0}(\cL_+(t,s),1;0)>0$, and $\cL_-(t,s)-\mL_{u,0}(\cL_-(t,s),1;0)<0$. For this, it is enough to verify that in the particular case of $\eps=\delta=0$ and $b=1$
\begin{eqnarray}\label{l1}
\hspace{-0.5cm}\cL_+\!\!&\!\!-\!\!&\!\!\Um\left(\!T_{b,\cL_+,\eps}\left(\!\sq{\Um\left(\!T_{b,\cL_+,\eps}\left(\!-\sq{{ \eps \cL_+\over 1+\eps}+{ \cL_+^2\over (1+\eps)^2}+{\cL_+ \delta \over 1+\eps}}{b}\!\right)\!;t,s\!\right)}{b}\! ;t,s\right) \!\right) \! >  0,\\ 
\label{l2}\hspace{-0.5cm} \cL_-\!\!&\!\!-\!\!&\!\!\um\left(\!T_{b,\cL_-,\eps}\left(\!\sq{\um\left(\!T_{b,\cL_-,\eps}\left(\!-\sq{{ \eps \cL_-\over 1+\eps}+{ \cL_-^2\over (1+\eps)^2}-{\cL_- \delta \over 1+\eps}}{b}\!\right)\!;t,s\!\right)}{b} \!;t,s\!\right) \!\right)  \! < 0
\end{eqnarray}
for all  $(t,s) \in \tilde{\cS}$ (we have omitted the arguments of functions $\cL_\pm$ above to make the notation less cumbersome). Since the left hand sides of the strict inequalities $(\ref{l1})$ and $(\ref{l2})$ are clearly continuous in $\eps$, $\delta$ and $b$, the same is true for sufficiently small $\eps$, $\delta$ and $\cB_- \le b \le \cB_+$. 

Inequalities $(\ref{l1})$ and $(\ref{l2})$ have been verified on a computer (see \cite{Progs}).

Next, we check that $\cB_+ > \mB_{u,\tau}(\lambda,\cB_+;e,\eps)$ and $\cB_- < \mB_{u,\tau}(\lambda,\cB_-;e,\eps)$. To verify $\cB_- < \mB_{u,\tau}(\lambda,\cB_-;e,\eps)$ we notice that
\begin{eqnarray}
\nonumber \hspace{-1.8mm}\mB_{u,\tau}(\lambda,\cB_-;e,\eps)\!\! &\!\!=\!\!&\!\!u \!\left( \! T_{\cB_-,\lambda,\eps} \! \left(-\sq{{\lambda \over 1+\eps} \left( \cB_- -e^2+\eps u(T_{\cB_-,\lambda,\eps}(e))- \tau \left(u(T_{\cB_-,\lambda,\eps}(e))\right) \right)  }{\cB_-}  \right)\! \right) \\
\nonumber \!\! & \!\! \ge \!\!&\!\!u \left(  T_{1+\sigma \eps^4,\lambda,\eps} \left(  -\sqrt{ 1+\sigma \eps^4-{\lambda \over 1+\eps}  \left( 1+\sigma \eps^4-\left( - \eps s {\alpha \over 2}   \right)^2 \right. }\right.\right.\\
\nonumber && \phantom{aaaaaaaaaaaa}\left.\left.\overline{\phantom{iia}+\left.\eps u  \left( T_{1+\sigma \eps^4,\lambda,\eps}\! \left(- \eps t {\alpha \over 2} \right) \right)- \tau \!\left( \!u \! \left( T_{1+\sigma \eps^4,\lambda,\eps}\! \left(e \right) \right) \right)  \right)    } \right) \right)\\
\nonumber  \!\!&\!\!\equiv\!\!& \!\!u(G(\lambda,s,e;\eps))\\
\nonumber \!\!&\!\!=\!\!&\!\!1+u'(G(\lambda,s,e;0)) \partial_\eps G(\lambda,s,e;0) \eps +{1 \over 2} \left[ u'(G(\lambda,s,e;0)) \partial^2_\eps G(\lambda,s,e;0) +\right.\\
\nonumber \!\!&\!\!\phantom{=}\!\!&\left.\hspace{1.5mm} + u''(G(\lambda,s,e;0)) \left(\partial_\eps G(\lambda,s,e;0)\right)^2 \right] \eps^2+O(\eps^3)\\
\nonumber \!\!&\!\!=\!\!&\!\! 1+t \partial_\eps G(\lambda,s,e;0) \eps \!+\!{1 \over 2} \!\!\left[ t \partial^2_\eps G(\lambda,s,e;0)\! +\!u''\!\!\left(\!-{1\over 2} \!\right) \!\left(\partial_\eps G(\lambda,s,e;0)\right)^2 \right] \eps^2\!+\!O(\eps^3).
 \end{eqnarray}

A straightforward but rather cumbersome calculation shows that $\partial_\eps G(\lambda,s,e;0)=O(\tau,\tau')$, and therefore $\partial_\eps G(\lambda,s,e;0) \equiv 0$ at $\tau=0$. At the same time, for $\tau=0$
$$\partial^2_\eps G(\lambda,s,e;0)=-\lambda {2 t^2-s^2 \over 32 (\lambda-1)^2}$$
which is positive for all $\cL_-(t,s) \le \lambda \le \cL_+(t,s)$ and $(t,s) \in \tilde{\cS}$. Therefore 
\begin{equation}\label{B_ineq_m}
 \mB_{u,0}(\lambda,\cB_-;e,\eps) \ge  1+{1 \over 2} t |\lambda| {2 t^2-s^2 \over 32 (\lambda-1)^2} \eps^2+O(\eps^3).
\end{equation}

For sufficiently small $\eps$ the right hand side of $(\ref{B_ineq_m})$ is strictly larger then $\cB_-=1+\sigma \eps^4$. This, together with the fact that $\mB_{u,\tau}(\lambda,\cB_-;e,\eps)$ is continuous in $\tau$ implies that the inequality $\cB_- < \mB_{u,\tau}(\lambda,\cB_-;e,\eps)$ holds for all sufficiently small $\tau$ and $\eps$.
 
To verify $\cB_+ > \mB_{u,\tau}(\lambda,\cB_+;e,\eps)$ we proceed in a similar way
\begin{eqnarray}
\nonumber \hspace{-1.8mm} \mB_{u,\tau}(\lambda,\cB_+;e,\eps)\!\! &\!\!=\!\!&\!\!u \!\left( \! T_{\cB_+,\lambda,\eps} \! \left(-\sq{{\lambda \over 1+\eps} \left( \cB_+ -e^2+\eps u(T_{\cB_+,\lambda,\eps}(e))- \tau \left(u(T_{\cB_+,\lambda,\eps}(e))\right) \right)  }{\cB_+}  \right)\! \right) \\
\nonumber \!\! & \!\! \le \!\!&\!\!u \left(  T_{1+C \eps^2,\lambda,\eps} \left(  -\sqrt{ 1+C \eps^2-{\lambda \over 1+\eps}  \left( 1+C \eps^2-\left( - \eps t {\alpha \over 2}   \right)^2 \right. }\right.\right.\\
\nonumber && \phantom{aaaaaaaaaaaa}\left.\left.\overline{\phantom{iia}+\left.\eps u  \left( T_{1+C \eps^2,\lambda,\eps}\! \left(- \eps s {\alpha \over 2} \right) \right)- \tau \!\left( \!u \! \left( T_{1+C \eps^2,\lambda,\eps}\! \left(e \right) \right) \right)  \right)    } \right) \right)\\
\nonumber  \!\!&\!\!\equiv\!\!&\!\! u(F(\lambda,t,e;\eps))\\
\nonumber \!\!&\!\!=\!\!&\!\!1+u'(F(\lambda,t,e;0)) \partial_\eps F(\lambda,t,e;0) \eps +{1 \over 2} \left[ u'(F(\lambda,t,e;0)) \partial^2_\eps F(\lambda,t,e;0) +\right.\\
\nonumber \!\!&\!\!\phantom{=}\!\!&\left.\hspace{1.5mm} + u''(F(\lambda,t,e;0)) \left(\partial_\eps F(\lambda,t,e;0)\right)^2 \right] \eps^2+O(\eps^3)\\
\nonumber \!\!&\!\!=\!\!&\!\! 1+t \partial_\eps F(\lambda,t,e;0) \eps \!+\!{1 \over 2} \!\left[ t \partial^2_\eps F(\lambda,t,e;0) \!+\!u''\left(\!-{1\over 2} \! \right) \left(\partial_\eps F(\lambda,s,e;0)\right)^2 \right] \eps^2 \!+\!O(\eps^3).
 \end{eqnarray}

Again, for $\tau=0$,
$$\partial_\eps F(\lambda,t,e;0)=0, \quad  \partial^2_\eps F(\lambda,t,e;0)=-\lambda{ 16 C (1-\lambda) +2 ts - t^2 \over 32 (1-\lambda)^2 },$$   
which is  positive for all $\cL_-(t,s) \le \lambda \le \cL_+(t,s)$ and $(t,s) \in \tilde{\cS}$ (where $s>t$), therefore
\begin{equation}\label{B_ineq_p}
\mB_{u,0}(\lambda,\cB_+;e,\eps) \le  1+{1 \over 2} t |\lambda|{ 16 C (1-\lambda) + 2 t s - t^2 \over 32 (1-\lambda)^2 } \eps^2+O(\eps^3) < 1+C \eps^2,
\end{equation}
if 
\begin{equation}\label{C2}
 {1 \over 2} t |\lambda|{ 16 C (1-\lambda ) +2 t s -t^2 \over 32 (1-\lambda)^2 }<C.
\end{equation}

Notice, that
$$\left[C-{1 \over 2} t |\lambda|{16 C (1-\lambda) + 2 t s - t^2 \over 32 (1-\lambda)^2 }\right]_{C=\left( t {\alpha(1,\lambda,0) \over 2 }\right)^2} = {t^2 (2-2 \lambda+\lambda s  )\over 32 (1-\lambda)^2}>0.$$

Therefore, conditions $(\ref{C1})$ and $(\ref{C2})$ are satisfied for all  $\cL_-(t,s) \le \lambda \le \cL_+(t,s)$ and $(t,s) \in \tilde{\cS}$ by some $C$ smaller, but sufficiently  close to $(t \ \alpha(b,\lambda,\epsilon)/2)^2$.

The solution $b$ is contained in the interval $\left( \mB_{u,\tau}(\lambda,\cB_-;e,\eps), \mB_{u,\tau}(\lambda,\cB_+;e,\eps) \right)$ which for sufficiently small $\eps$ is a subset of  
$$(\hat{\cB}_-,\hat{\cB}_+) \equiv \left(1+{1 \over 2} t |\lambda| {2 t^2-s^2 \over 32 (\lambda-1)^2} \eps^2, 1+ {1 \over 2} t |\lambda|{16 C (1-\lambda) + 2 t s-t^2  \over 32 (1-\lambda)^2 } \eps^2 \right).$$

Notice, that for $C=0$
$$\hat{\cB}_+ - \hat{\cB}_-=-\lambda t {16 C (1-\lambda) +2 ts +s^2-3 t^2  \over 64 (1-\lambda)^2}$$ 
which is positive for all $(t,s) \in \tilde{\cS}$ where $s>t$. Therefore the interval $(\hat{\cB}_-,\hat{\cB}_+)$ is non-empty.

$\Box$

\subsection{Proof of part ii) of Proposition  $\ref{central_prop}$.}\label{part2}

Differentiate $\oT[u]$ with respect to $x$:
\begin{equation}
\nonumber \oT[u]'(x)= {\alpha^2 \over 1+\eps}u' \left( \Tble \left( \Veut( \Tble^{-1}(x)) \right)\right) {u'\left(\Tble \left(-\sqrt{w(x) }\right)\right) \over  4 \Veut(\Tble^{-1}(x) )} {w'(x) \over \sqrt{w(x)}  },
\end{equation}
where  $w$ is the function defined in $(\ref{w_function})$. On the real line 
$$ \mw \le  w \le \mW \quad {\rm and} \quad \mv \le \Veut \circ \Tble^{-1} \le \mV,$$
where
\begin{eqnarray}
\nonumber \mW (x;t,s) &=&{b- {\lambda \over 1+\eps} \left(b-\Tble^{-1}(x)^2+ \eps \Um(x;t,s) +\delta \right)}, \\
\nonumber \mw (x;t,s) &=&{b- {\lambda \over 1+\eps} \left(b-\Tble^{-1}(x)^2+ \eps \um(x;t,s) -\delta \right)},\\
\nonumber \mV (x;ts) &=& \sq{u  \left( \Tble  \left( - \sqrt{ {\mw(x;t,s)}}\right)\right)}{b}, \\
\nonumber  \mv (x;ts) &=& \sq{u  \left( \Tble  \left( - \sqrt{\mW(x;t,s)}\right)\right)}{b}
\end{eqnarray}
are upper and lower bounds on the corresponding functions. Notice that 
$$ u'(x) \le \Theta_2'(\mF_2(\Phi_1(x);t,s) ) \mDf(\Phi_1(x);t,s) \Phi_1'(x) \equiv \mDu(x;t,s)
$$
where
$$ \mDf(x;t,s) \equiv \eta(x-c_1) \mF_2(x,t,s) {(1-c_1)  \over  (x-c_1) (1-x) } + \eta(c_1-x)  \mF_2(x,t,s) {(1+c_1)  \over  (x-c_1) (1+x) }$$
is an upper bound on derivatives on $\cA_1({\bf c})$ that follows from $(\ref{first_der})$ ($\eta$ is the Heaviside function). Therefore,
$$\oT[u]'(x) \le  {\alpha^2 \over 1+\eps} \mDu \left( \Tble \left(\mv(x;t,s) \right)\right) {\mDu \left(\Tble \!  \left(\! -\sqrt{\mW(x;t,s) } \right) \right)\!\! \over  4 \mv(x;t,s)} \cdot {\!-2 \alpha^{-1}\Tble^{-1}(x) -\varepsilon (\omega+\rho x) \!\! \over \sqrt{\mw(x;ts)}  }.$$
We finally verify on the computer (see \cite{Progs}) that the right hand side of the above inequality is less than $\omega+\rho x$ for all $x \in (0,r_1)$ and sufficiently small $\eps$ and $\varepsilon$.

$\Box$

\subsection{Proof of part iii) and iv) of Proposition  $\ref{central_prop}$.}\label{part23}

Suppose that $\theta \mapsto \partial \cD(\theta)$ and $\theta \mapsto \partial \cE(\theta)$ are some convenient parametrization of the boundaries, such that $\partial \cD \cap \fC_+$ is parametrized by $\theta \in (0,\pi)$, while $\partial \cD \cap \fC_-$ is parametrized by $\theta \in (-\pi,0)$, and similarly for $\cE$. Let,  again, $w$ be the function defined in $(\ref{w_function})$. Denote $H$ the preimage of the ray $\left( w(\Tble(e)),+\infty \right)$ in $\Tble^{-1}(\cD)$.  

First, we would like to find  a bound on  $W(\theta) \equiv  T_{b,\lambda,\eps} \left(-\sqrt{w(  \partial \cD(\theta) )   }\right)$  for $0 \le \theta \le \pi$. For every fixed $\theta$,  $W(\theta)$ is contained in the set $\cW(\theta)$  bounded by the  curves 
\begin{eqnarray}
\nonumber W_\cE(\theta,p)&=& T_{b,\lambda,\eps} \left(-\sq{{\lambda \over 1+\eps} \left(b-T^{-1}_{b,\lambda,\eps}(\partial \cD (\theta) )^2 +\varepsilon \partial \cE(p)-\tau(\partial \cE(p)) \right)}{b}\right), \quad 0 \le p \le \pi \\
\nonumber W_\Re (\theta,p)&=& T_{b,\lambda,\eps} \left(-s_b\left({\lambda \over 1+\eps} \left(b-T^{-1}_{b,\lambda,\eps}(\partial \cD (\theta) )^2 +\varepsilon (p \partial \cE(0) +(1-p) \partial \cE(\pi) ) \right.  \right. \right. \\
\nonumber & & \phantom{aaaaaaaaaaaaaaa} \left.  \left.\phantom{ 1\over \eps}+ \tau( p \partial \cE(0) +(1-p) \partial \cE(\pi)  )\right)  \right), \quad 0 \le p \le 1.
\end{eqnarray}

\begin{figure}[t]
 \begin{center}
\begin{tabular}{c c}
\resizebox{70mm}{!}{\includegraphics[angle=-90]{fig4-1.eps} } & \resizebox{70mm}{!}{\includegraphics[angle=-90]{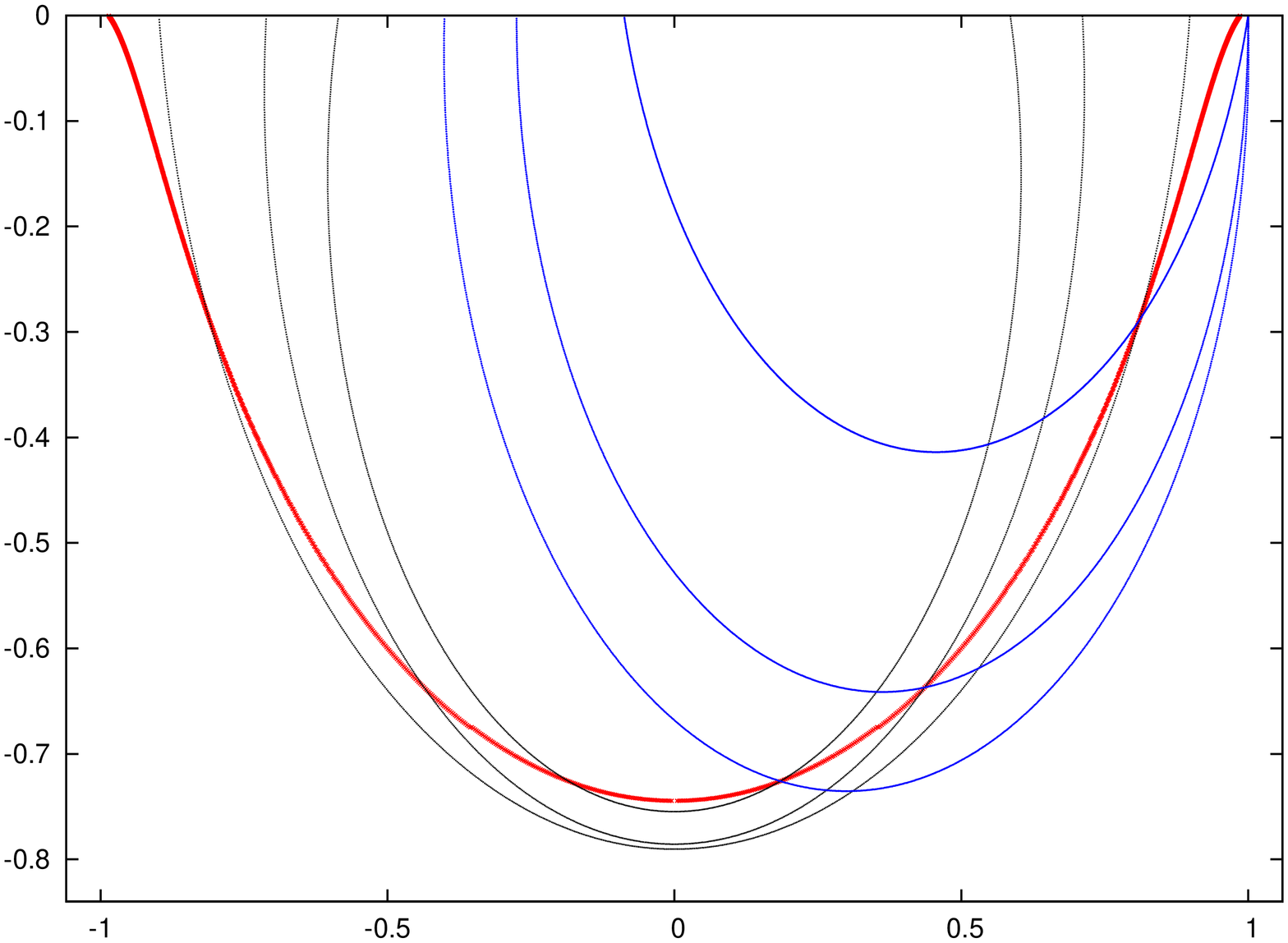}} \\
a) & b) 
\end{tabular}
\caption{a) Orbit of the set $\Phi_1(\cW(\theta))$ for  $0 \le \theta \le \pi$ (red). An example of a cover given for six boundary points: each point of the orbit is in the intersection of two Poincar\'e neighborhoods (black and blue) in collections $\cup_n \cD_+(J^+_n,\theta^+_n)$ (black and blue lines in the upper half plane) and  $\cup_m \cD_-(J^-_m,\theta^-_m) $ (black and blue lines in the lower half plane); b) The boundary of the set $\cN$ is given in red. An example of a cover for three boundary points:  each point of the boundary is in the intersection of two Poincar\'e neighborhoods (black and blue) in the collection $\cH$.}  
\label{Set_W}
\end{center}
\end{figure}

As before, we consider the case $\varepsilon=\tau=0$, by continuity of all involved functions, the claim will also hold for sufficiently small $\varepsilon$ and $\tau$. Recall, that  $u= \Theta_2 \circ f_2 \circ \Phi_1=\Theta_3 \circ f_3 \circ \Phi_1$. We first cover the set $\Phi_1(\cW(\theta))$, $0 \le \theta \le \pi$ by a collection of  Poincar\'e half-neighbourhoods 
$$\cP=\left(\cup_n \cD_+(J^+_n,\theta^+_n) \right) \cup \left(\cup_m \cD_-(J^-_m,\theta^-_m)\right)$$
for some appropriately chosen $J^+_n=(l^+_n,r^+-_n)$, $J^-_m=(l^-_m,r^-_m)$ and $\theta^-_n$, $\theta^+_m$ (cf. \ref{Set_W}), then according to Lemma $\ref{Epstein_lemma}$, the sets $f_k(\Phi_1(\cW(\theta)))$, $k=2,3$, $0 \le \theta \le \pi$, are contained in 
$$\cU_k(t,s)=\left( \cup_n \cD_+(\tilde{J}^+_{n,k},\theta^+_n) \right) \cup \left(\cup_m \cD_-(\tilde{J}^-_{m,k},\theta^-_m) \right), \quad k=2,3$$
where  $\tilde{J}^\pm_{i,k}=(\mf_k(l^\pm_i;t,s),\mF_k(r^\pm_i;t,s))$, $k=2,3$. Set $\cV(t,s) \equiv \Theta_2(\cU_2(t,s)) \cap \Theta_3(\cU_3(t,s))$. 

The choice of neighborhoods $\cD_\pm(J^\pm_n,\theta^\pm_n)$ is implemented on a computer via an automatized procedure (see  \cite{Progs}): the neighborhoods are constructed so that every point $z$ of the curve $\Phi_1(\cW(\theta))$ lies in the intersection of two such neighborhoods  $\cD_\pm(J^\pm_{n'},\theta^\pm_{n'})$ and  $\cD_\pm(J^\pm_{n''},\theta^\pm_{n''})$, then $f_k(z)$ lies in the intersection of  $\cD_\pm(\tilde{J}^\pm_{n'},\theta^\pm_{n'})$ and  $\cD_\pm(\tilde{J}^\pm_{n''},\theta^\pm_{n''})$.

 We next construct the set 
$$\cM(b,\lambda;t,s)=-{\rm sign}\left(\Im\left(b-\cV(t,s) \right)  \right)\sqrt{ b-\cV(t,s)},$$
which is a bound on $\Veut(\Tble^{-1}(\cD))$, and verify that it is contained in $\Tble^{-1}(\cD)$. Similarly to $V_u$ (cf Prop. $\ref{central_prop_F}$, part $ii)$), $\Veut$ is continuous across $H$ and holomorphic in $\Tble^{-1}(\cD)  \setminus H$; by Morera's theorem  it is holomorphic in $\Tble^{-1}(\cD)$.

Next, construct
$$\cN(b,\lambda;t,s)=\Phi_1\left( \Tble \left(-\cM(b,\lambda);t,s \right)\right)$$
and cover it with another pair of collections of Poincar\'e  half-neighbourhoods 
$$\cH=\left(\cup_n \cD_+(I^+_n,\phi^+_n) \right) \cup \left(\cup_m \cD_-(I^-_m,\phi^-_m)\right).$$
Set
$$\tilde{\cH}_k(t,s)=  \left( \cup_n \cD_+(\tilde{I}^+_{n,k},\phi^+_n) \right) \cup \left(\cup_m \cD_-(\tilde{I}^-_{m,k},\phi^-_m) \right), \quad k=2,3, $$
where $\tilde{I}^\pm_{i,k}=(\mf_k(l^\pm_i;t,s),\mF_k(r^\pm_i;t,s))$, $k=2,3$. Finally, the set
$$\cX(b,\lambda;t,s)=\lambda^{-1} \Theta_2 \left(  \tilde{\cH}_2(t,s)\right) \cap  \lambda^{-1} \Theta_3 \left(  \tilde{\cH}_3(t,s) \right),$$
which is  a bound on $\oT[u](\cD)$,  is verified to be contained  in $\cE$. This shows that $\cT_{\eps,\tau}[u]$ is in $\cA(\cD,\cE,{\bf \mc})$ whenever $u \in \cA(\cD,\cE,{\bf \mc})$.

$\Box$

\section{Appendix A: New a-priori bounds on ${\fR}$}\label{new_bounds}

In this subsection we will use {\it a-priori} bounds on ${\cA}_1({\bf c}))$ to produce  better bounds on a subset of functions bounded on $(-1,1)$ by a constant.

As before, we denote $(t,s)=(u'(-1/2),u'(0) )$ for a function $u \in  \cA(\cD,\cE, {\bf \mc})$, ${\bf \mc}=(-1/2,0,0,1)$. Recall that $u= \Theta_k \circ f_k \circ \Phi_1$, $k=2,3$, where $f_k \in \cA_1({\bf c^k})$, ${\bf c^k}=(\Phi_1(c_1),\Phi_1(c_2),\Phi_k(c_3),\Phi_k(c_4))$ (note, we will be using the superscript $k$ on functions and numbers, whenever convenient, to avoid double subscripts, these by no means signify raising to a power). Therefore, the following are the derivatives of $f_k$ at points $c_1$ and $c_2$: 
$$T_k(t)={t \over \Theta'_k(c^k_3) \Phi'_1(-1/2)}, \quad S_k(s)={s \over \Theta'_k(c^k_4) \Phi'_1(0)}.$$

Now, recall that $f'(x)$ is convex, and therefore, using $(\ref{second_der})$,
\begin{eqnarray}\label{2der_bounds}
 & &\min_{x \in [c_1,c_2]} f''_k(x)  \ge -2 {(c_2-x) T_k(t)+(x-c_1) S_k(s) \over (c_2-c_1) (1+c_1)} \equiv m_k(x,t,s), \\
 & &\max_{x \in [c_1,c_2]} f''_k(x) \le 2 { (c_2-x) T_k(t) +(x-c_1) S_k(s) \over (c_2-c_1) (1-c_2)} \equiv M_k(x,t,s).
\end{eqnarray} 
Now, fix $t$ and $s$, and consider the function $y_k(x)=T_k(t)+\int^x_{c_1} m_k(z,t,s) d z$. Suppose, that the line $w_k(x)=S_k(s)+n_k(t,s) (x-c_2)$ intersects $(x,y_k(x))$ at point $x_k(t,s)$, and $n_k(t,s)$ is such that the following holds:
$$c_4^k-c_3^k=\int^{c_2}_{c_1} \my_k(z) d z, \quad 
\my_k(x)=\left\{y_k(x), \quad c_1 \le x \le x_k(t,s) \atop w_k(x), \quad x_k(t,s) \le x \le c_2 \right. .$$
First, notice, that any curve $(x,f'_k(x))$ on $(c_1,c_2)$ with end points $(c_1,t)$ and $(c_2,s)$ can not intersect $(x,y_k(x))$, and has to intersect $(x,w_k(x))$ somewhere on $(x_k(t,s),c_2)$ once ($f'_k(x)$ is convex), for if it does not then $\int_{c_1}^{c_2} f'_k(z) d z  \ne c_4-c_3$. It is also clear that 
\begin{equation}
f_k(x)  \ge c^k_3+ \int_{c_1}^x y_k(z) d z \equiv f^k_2(x;t,s), \quad  x \in \left[c_1,c_2 \right].
\end{equation}

One can repeat a similar argument for $Y_k(x)=S_k(s)+\int^x_{c_2} M_k(z,t,s) d z$ and $W_k(x)=T_k(t)+N_k(t,s)(x-c_1)$ that intersect at $X_k(t,s)$ to get 
$$f_k(x) \le c^k_4-\int_x^{c_2} \mY_k(z) d z \equiv F^k_2(x;t,s),  \quad  x \in \left[c_1,c_2 \right], \quad  \mY_k(x)=\left\{Y_k(x), \quad  X_k(t,s) \le x \le c_2 \atop W_k(x), \quad c_1 \le x \le X_k(t,s) \right. .$$

To obtain an  upper bound on $(-1,c_2)$ and a lower bound on $(c_2,1)$, we recall that the positivity of the Schwarzian derivative for functions in ${\cA}_1({\bf c^k}))$ together with the positivity of all $f^{(n)}_k$ for odd $n$ implies that for all $x \in (-1,1)$ 
\begin{equation}\label{S_der}
f'''_k(x) \ge {3 f''_k(x)^2 \over 2 f'_k(x)},
\end{equation}
and consequently,
$$f''_k(x) \le f''_k(c_1)+{3 \over 2} \int_{c_1}^x {f''_k(y)^2 \over  f'_k(y)} d y,$$
for all $x \in (-1,c_1)$, the equality being realized by the the solution 
$$f'_k(x)={4 f'_k(c_1)^3 \over ( -f''_k(c_1) (x-c_1) +2 f'_k(c_1) )^2}$$
of equation $(\ref{S_der})$. Therefore, 
$$f_k(x) \le \int_{c_1}^x {4 T_k(t)^3 \over ( -f''_k(c_1) (x-c_1) +2 T_k(t) )^2}, $$
for all $x \in (-1,c_1)$, the maximum of the right hand side being realized by the maximum admissible $f''_k(c_1)$ which can be obtained from the condition 
\begin{equation}\label{S_derr}
 {4 T_k(t)^3 \over ( -f''_k(c_1) (c_2-c_1) +2 T_k(t) )^2} = S_k(s).
\end{equation}
We denote $Z_k(t,s)$ the solution $f''_k(c_1)$ of this equation, then 
$$f_k(x) \le c_3^k+{4 T_k(t)^3 \over Z_k(t,s) }\left({1 \over 2 T_k(t) +Z_k(t,s) (c_1-x) }-{1 \over 2 T_k(t) } \right) \equiv F^k_1(x;t,s), \quad x \in (-1,c_1).$$ 

In a similar way
$$  f_k(x)  \ge c_4^k+{4 S_k(t)^3 \over X_k(t,s) }\left({1 \over 2 S_k(t) +X_k(t,s) (c_2-x) } - {1 \over 2 S_k(t) } \right)=f^k_3(x;t,s), \quad x \in (c_2,1),$$
here $X_k(t,s)$ solves
\begin{equation}\label{S_derrr}
 {4 S_k(t)^3 \over ( -X_k(t,s) (c_1-c_2) +2 S_k(t) )^2} = T_k(s).
\end{equation}

Finally, suppose that $ \mm_k \le f_k(x) \le \mM_k$ on the real slice of its domain (this is certainly true if $f_k  \in  \cA_1({\bf c^k})$). Consider the line  $(x,S_k(s)+\mK_k (x-c_2))$ where $\mK_k$ is such that 
$$\int_{c_2}^1  S_k(s)+\mK_k (x-c_2) d x =\mM_k-c_4^k,$$
that is 
$$\mK_k=2{\mM_k-c_4^k \over (1-c_2)^2}-{S_k(s) \over 1-c_2}.$$ 
Since $f'_k(x)$ is convex, the curve $(x,f'_k(x))$ intersects the line $(x,S_k(s)+\mK_k (x-c_2))$ strictly once on $(c_2,1)$. Convexity of $f'_k(x)$ implies that 
$$\int_{c_2}^x f'_k(y) dy<\int_{c_2}^x  S_k(s)+\mK_k (y-c_2) d y, \quad x \in (c_2,1),$$
that is
\begin{equation}
f_k(x) \le c^k_4+S_k(s)(x-c_2) +(\mM_k-c_4^k-S_k(s) (1-c_2) ) {(x-c_2)^2 \over (1-c_2)^2 } \equiv F^k_3(x;t,s), \quad x \in (c_2,1).
\end{equation}

A similar argument on $(-1,c_1)$ demonstrates that 
\begin{equation}
f_k(x)\ge c^k_3 -T_k(t) (c_1-x) + (T_k(t) (1+c_1)+\mm_k-c_3^k){(x-c_1)^2 \over (1+c_1)^2}  \equiv f^k_1(x;t,s), \quad x \in (-1,c_1).
\end{equation}

Finally, $\mf_k(x;t,s) \le f_k(x) \le \mF_k(x;t,s)$ on $(-1,1)$, where
\begin{equation}\label{f_bounds}
\mf_k(x;t,s)=\!\left\{\!\begin{array}{cc} f^k_1(x;t,s),  x \in & \left(-1,c_1 \right) \\  f^k_2(x;t,s)),  x \in & \left(c_1,c_2 \right) \\  f^k_3(x;t,s)),  x \in & \left(c_2,1\right) \end{array}  \!\right., \quad    \mF_k(x;t,s)=\!\left\{\!\begin{array}{cc} F^k_1(x;t,s),  x \in & \left(-1,c_1\right) \\  F^k_2(x;t,s)),  x \in & \left(c_1,c_2 \right) \\  F^k_3(x;t,s)),  x \in & \left(c_2,1\right) \end{array} \! \right.
\end{equation}

Bounds $(\ref{f_bounds})$ transferred to the space $\cA(\cD,\cE;{\bf \mc})$ will be denoted $\um$ and $\Um$:
\begin{eqnarray}
\label{uu_bounds} \um(x;t,s) &\equiv& \max \left(\Theta_2(\mf_2 (\Phi_1(x);t,s)),\Theta_3(\mf_3 (\Phi_1(x);t,s))\right),\\
 \Um(x;t,s) &\equiv& \min \left( \Theta_2(\mF_2 (\Phi_1(x);t,s)),  \Theta_3(\mF_3 (\Phi_1(x);t,s)) \right).
\end{eqnarray}

\section{Appendix B: Set of realizable $(u'(-1/2),u'(0))$}\label{ts}

In this subsection we will describe the set $\cS$ of realizable $t=u'(-1/2)$ and $s=u'(0)$ whenever $u \in \cA(\cD, \cE, {\bf  c})$, and its subset $\tilde{\cS} \subset \cS$ invariant under $\oT$.

 Write  $u= \Theta_k \circ f_k \circ \Phi_1$, $k=2,3$,  $f_k \in \cA_1({\bf c})$ as before. Since $f_k(x) \le F^k_1(x;t,s)$ on $(-1,c_1)$ (see Subsection $\ref{new_bounds}$) we have $ -1 \le  F^k_1(-1;t,s)$. The relevant (positive) solution $s=s(t)$ of  this equation  will be denoted  by $\cZ_k(t)$. Similarly, $f^k_3(1;t,s)  \le 1$.  The relevant solution $s=s(t)$ will be denoted by $\cC_k(t)$. 
 We have obtained symbolic (and not just numeric) expressions for $\cZ_k(t)$ and $\cC_k(t)$  using the Maple software package. The set bounded by these curves is the set $\cS$ of admissible values $(t,s)$.

\begin{figure}[t] 
 \begin{center}
 \resizebox{70mm}{!}{\includegraphics[angle=-90]{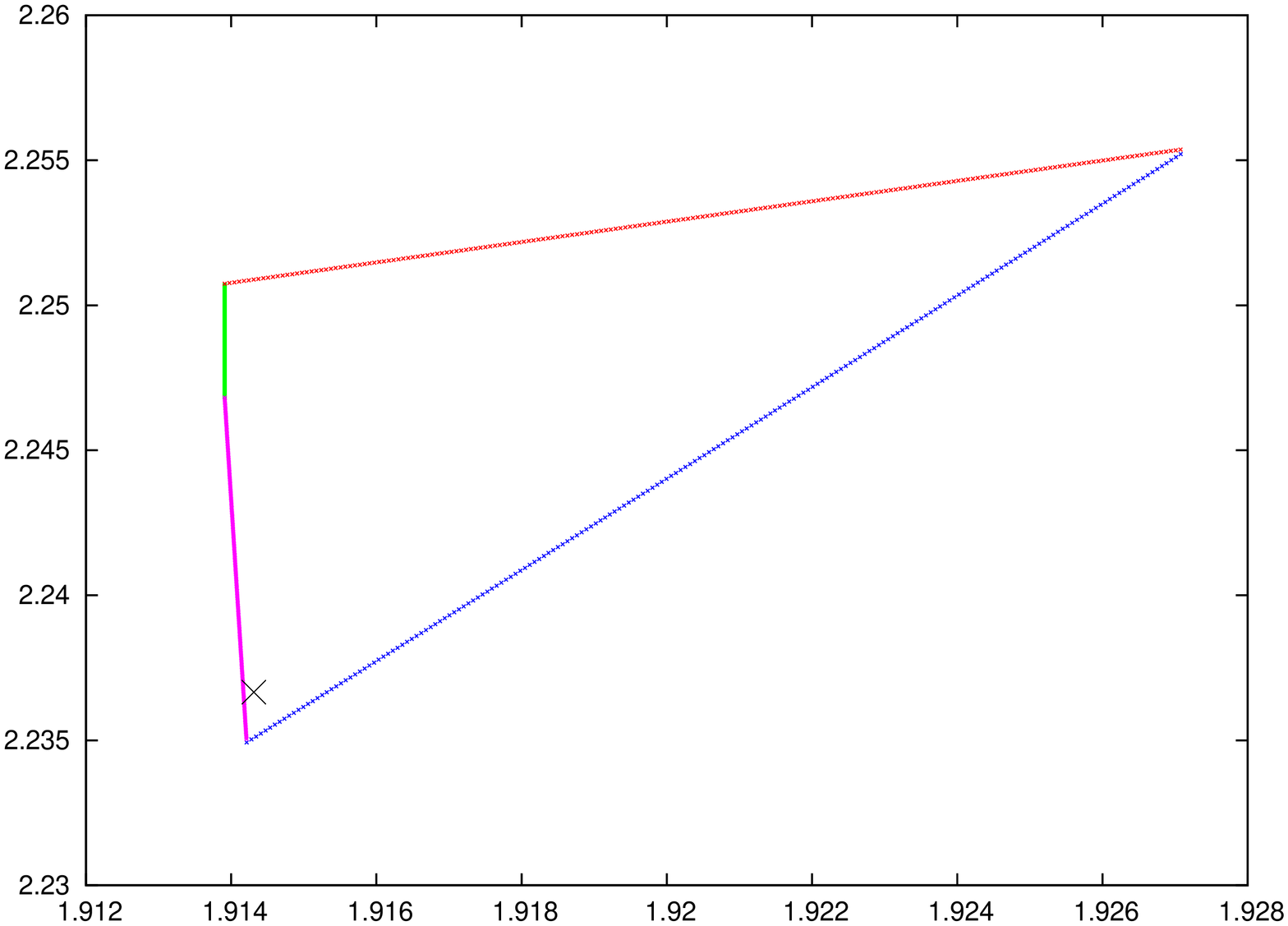}}  
\caption{ Invariant set $\tilde{\cS}$ bounded by curves $\cZ_2$ (blue), $\cZ_3$ (magenta), $t=t^*-0.0004$ (green) and $\cC_2$ (red). The cross marks the location of $(t^*,s^*)$. \label{set_S}}
\end{center}
\end{figure}


We can further restrict the set of admissible $(t,s)$ if we notice that 
$$\cT[u]'(-1/2) ={- \alpha(1,\lambda,0) t s  \over 2 \lambda \beta(1,\lambda,0) } = -{t s  \over 4 \lambda \beta(1,\lambda,0)^2 }= {t s  \over 4 \lambda (\lambda-1)} \equiv  \mT(\lambda,t,s).$$

Denote $\tilde{\cS}$ the subset of $\cS$  to the left of the line $t=t^*-0.0004$.
 We have  verified on the computer that for all $(t,s) \in \tilde{\cS}$, all $\cL_-(t,s) \le  \lambda \le \cL_+(t,s)$, $b=1$, $\eps=0$, 
\begin{equation}\label{S_inv}
\mT(\lambda,t,s) > t^*-0.0004,
\end{equation} 
(see \cite{Progs}).

We have shown in Prop. $\ref{central_prop}$, part $iv)$, that if the derivatives $(t,s)$ for a function  $u \in \cA(\cD, \cE, {\bf  c})$ are in $\tilde{\cS}$, then $\oT[u] \in  \cA(\cD, \cE, {\bf  c})$, that is  $\left( \cT_{\eps,\tau}[u]'(-1/2),\cT_{\eps,\tau}[u]'(0)  \right) \in \cS$. This, together with the strict inequality $(\ref{S_inv})$, implies that the subset $\tilde{\cS}$ is invariant under the map $(t,s) \mapsto \left( \cT_{\eps,\tau}[u]'(-1/2),  \cT_{\eps,\tau}[u]'(0)  \right)$ for nonzero $\eps$ and $\tau$.  This subset is depicted in Fig. \ref{set_S}.


\section{Acknowledgements}

The author would like to cordially thank  Michael Benedicks (KTH, Stockholm) and Hans Koch (University of Texas, Austin, Texas) for many useful discussions on the subject. 

\singlespacing

\end{document}